\UseRawInputEncoding
\documentclass[A4paper]{article}
\usepackage{amsfonts}
\usepackage{amsxtra}
\usepackage{amsthm}
\usepackage{amssymb}
\newtheorem{fed}{\textbf{Definition}}[section]
\newtheorem{thm}[fed]{\textbf{Theorem}}
\newtheorem{lemma}[fed]{\textbf{Lemma}}

\newtheorem{rem}[fed]{\textbf{Remark}}
\newtheorem{prop}[fed]{\textbf{Proposition}}
\newtheorem{cor}[fed]{\textbf{Corollary}}
\usepackage{amscd,amssymb,bbm,graphicx,epsfig,psfrag,epic,eepic,latexsym,amsmath}
\usepackage{amsmath}
\usepackage[arrow,matrix,curve]{xy}
\usepackage{mathrsfs}
\usepackage{graphicx}
\usepackage{color}
\usepackage{hyperref}

\newcommand{\N}{\mathbb{N}}

\newcommand{\Z}{\mathbb{Z}}
\newcommand{\T}{\mathbb{T}}
\newcommand{\R}{\mathbb{R}}
\newcommand{\C}{\mathbb{C}}

\begin{document}

\title{{Generalized Periodic Orbits of Time-Periodically Forced Kepler Problem Accumulating the Center and of Circular and Elliptic Restricted Three-Body Problems}}
\author{Lei Zhao\footnote{University of Augsburg, Augsburg, Germany: lei.zhao@math.uni-augsburg.de}}
\maketitle

\abstract 
In this paper, we consider a time-periodically forced Kepler problem in any dimensions, with an external force which we only assume to be regular in a neighborhood of the attractive center. We prove that there exist infinitely many periodic orbits in this system, with possible double collisions with the center regularized, which accumulate the attractive center. The result is obtained via a localization argument combined with a result on $C^{1}$-persistence of closed orbits by a local homotopy-streching argument. Consequently, by formulating the circular and elliptic restricted three-body problems of any dimensions as time-periodically forced Kepler problems, we obtain that there exists infinitely many periodic orbits, with possible double collisions with the primaries regularized, accumulating to each of the primaries. 

Appendix A coauthored with U. Frauenfelder.

\endabstract

\tableofcontents

\section{Introduction}
An essential part of perturbation theory in celestial mechanics is based on perturbations of uncoupled Kepler problems. The Kepler problem in $\R^{d}, d \ge 2$ is \emph{properly-degenerate} in the sense that it has $d$-degrees of freedom while all of its non-singular bounded orbits are closed and thus carries only one non-trivial frequency. Therefore to understand concrete problems as perturbations of Kepler problems typically one has to study higher order effects. 

On the other hand, while looking for periodic orbits of a fixed period, a different approach can be taken, as has been shown in \cite{Boscaggin--Ortega--Zhao}, which considered a periodically forced Kepler problem. By the third Kepler law, the period of periodic orbits of the Kepler problem changes with respect to the energy and this gives the normal non-degeneracy of the set of periodic orbits with fixed period in the extended phase space. This, combined with regularization techniques in extended phase space, makes it possible to apply some global methods from symplectic geometry to conclude the existence of periodic orbits after regularization, which correspond to periodic orbits in a generalized sense of the original problem.

In this paper, we consider the same type of model as in \cite{Boscaggin--Ortega--Zhao}, namely the time-periodically forced Kepler problem in $\R^{d}$ for $d \in \Z_{+}$:
\begin{equation}\label{eq: eq motion perturbed Kepler varepsilon=1}
\ddot{q}=-\dfrac{q}{|q|^{3}} +\nabla _{q} U(q, t),
\end{equation}
in which $q \in \R^{d}$. The function $U(q, t, \varepsilon)$ is a time-periodic function which we assume is $C^{\infty}$-regular in a neighborhood of the origin. Outside this neighborhood we do not make any further assumptions. Also, when $U(q, t, \varepsilon)$ is independent of $t$, we may fix any positive period. By normalization in the variables $(q, t)$ we may assume that $U(q, t, \varepsilon)$ is 1-periodic in $t$. 

In many concrete models the problem is described via an external parameter. Thus it is also helpful to consider the following interpolating system of the form
\begin{equation}\label{eq: eq motion perturbed Kepler}
\ddot{q}=-\dfrac{q}{|q|^{3}} +\varepsilon \nabla _{q} U(q, t, \varepsilon),
\end{equation}
in which $q \in \R^{d}$ and $U(q, t, \varepsilon)$ is a time-periodic function which is $C^{\infty}$-regular in a neighborhood of the origin depending continuously on a parameter $\varepsilon \in [0, 1]$ such that $U(q, t):=U(q, t, 1)$. 


This system is a non-autonomous Hamiltonian system with Hamiltonian
$$F_{\varepsilon}(p, q, t)=\dfrac{|p|^{2}}{2}-\dfrac{1}{|q|}+ \varepsilon U(q, t, \varepsilon).$$
Without loss of generality we may assume 
\begin{equation}\label{assumption: U}
U(0, t, \varepsilon)=0.
\end{equation} 
Indeed when this is not the case, it is enough to take $U(q,t, \varepsilon)-U(0,t, \varepsilon)$ in place of $U(q,t)$ in the above expression. As a consequence, $U(q, t, \varepsilon)$ admits a Taylor expansion with $t$-periodic coefficients of the form
$$U(q, t) = \sum_{i=1}^{d} c_{i}(t, \varepsilon) q_{i} +O(|q|^{2}), $$
and consequently in a neighborhood of the origin there holds
\begin{equation}\label{eq: estimate on U}
|U(q, t)| \le C_{0} |q| +O(|q|^{2}, t)
\end{equation}
for some constant $C_{0}>0$.

Generalized periodic solutions are solutions which become periodic after possible isolated double collisions with the attractive center being regularized. Precisely a generalized solution of Eq. \eqref{eq: eq motion perturbed Kepler} is a continuous function $\R \to \R^{d}, t \mapsto q(t)$ which satisfies Eq. \eqref{eq: eq motion perturbed Kepler} except for a discrete set $\mathcal{Z}:=\{t \in \R: q(t)=0\}$ of collision instants. Moreover, for any $t_{0} \in \mathcal{Z}$, the limits
$$\lim_{t \to t_{0}} \dfrac{u(t)}{|u(t)|}, \qquad \hbox{ and } \qquad \lim_{t \to t_{0}} \Bigl(\dfrac{1}{2} |\dot{q}(t)|^{2}-\dfrac{1}{|u(t)|}\Bigr) $$ 
of collision direction and collision energy exist. The last condition gives a condition for patching solutions which runs into and out of the collisions. 

The result in \cite{Boscaggin--Ortega--Zhao} concerns the number of generalized periodic orbits of Eq.\eqref{eq: eq motion perturbed Kepler}. With the remark from \cite{OrtegaZhao}, it is sufficient that the function $U(q, t, \varepsilon)$ be defined in a neighborhood of the origin for the $q$-variable as we now assume. The main result is:

For $d=2,3$, given any positive integer $l \in \N_{+}$, there exists some $\varepsilon_{*}>0$, such that for all $\varepsilon \in [0, \varepsilon_{*})$, the equation \eqref{eq: eq motion perturbed Kepler} has at least $l$ generalized periodic solutions. 

The result was established with regularization techniques and application of a theorem of Weinstein \cite{Weinstein} concerning bifurcations of periodic orbits from a non-degenerate periodic manifold when the perturbation is $C^{2}$-small in its neighborhood. The external small parameter in Eq. \eqref{eq: eq motion perturbed Kepler} was technically needed to control the $C^{2}$-norm of the perturbation. 

We mention two recent non-perturbative results concerning equations of the form \eqref{eq: eq motion perturbed Kepler varepsilon=1}:

For $d=2$  and $U(q, t)=O(|q|^{\alpha})$ for some $\alpha \in (0, 2)$ when $|q| \to \infty$, the existence of at least one generalized periodic orbit has been shown in \cite{Boscaggin-D'Ambrosio-Papini}. 

For $d=2,3$ and $U(q, t)=\langle p(t), q \rangle$, the existence of infinitely many generalized periodic orbits has been established in \cite{Verzini-Ortega-Barutello}.

These results were obtained via variational methods. Besides these, many results on the existence and multiplicity of periodic solutions with more general singular potentials has been obtained in e.g. \cite{Bahri-Rabinowitz}, \cite{Ambrosetti-Struwe} with variational methods. Some of these results have been collected in the book \cite{Ambrosetti-Coti Zelati} in which one may find many further references.

In this article we aim to establish the following theorem which uniformly treat the problem in all dimensions with general external force. In dimension 2 or 3, this enumeratively improves the result from \cite{Boscaggin--Ortega--Zhao} and allows more general forcing with more information about the space locations of the generalized periodic orbits to be found as compared to \cite{Verzini-Ortega-Barutello}:
\begin{thm}\label{thm: main} For all $d \in \Z_{+}$, there exists infinitely many generalized periodic orbits in the system \eqref{eq: eq motion perturbed Kepler varepsilon=1}, resp. \eqref{eq: eq motion perturbed Kepler}, accumulating the attractive center.
\end{thm}


The result is, in a sense, a typical result of the lunar regime: The Newtonian attraction is singular at the attractive center and therefore in a sufficiently small neighborhood eventually dominates the additional external regular force. Concretely, after embedding the system into the zero-energy hypersurface of an autonomous system and regularization, the generalized $1$-periodic orbits of Eq. \eqref{eq: eq motion perturbed Kepler} are transformed into certain closed orbits of the regularized system in extended phase space. When the external force is disregarded, these closed orbits form an infinite sequence of periodic manifolds $\Lambda_{n}$ for those with initial  prime period $1/n$ in the unperturbed system. These periodic manifolds are non-degenerate and are distinguished not only by their prime periods but also by their action values. We shall show that, when the external force is added, some closed orbits bifurcating from $\Lambda_{n}$ continue to exist, for all sufficiently large $n$. An estimate on their action values then shows that there are infinitely many of them. The persistence of closed orbits bifurcating from $\Lambda_{n}$ for $n$ large will be established via a local homotopy-streching argument from Rabinowitz-Floer homology theory adapted to our situation.  Indeed after a rescaling argument one may fix a periodic manifold and use the rescaling parameter to control the $C^{1}$-norm of the perturbations. In general this fails to bring a control for the $C^{2}$-norm of the perturbations and therefore in general Weinstein's theorem from \cite{Weinstein}, which was used in \cite{Boscaggin--Ortega--Zhao}, is not applicable in our setting. We consider this as an illustration of how recent developments of symplectic topology enrich the understanding of problems from classical and celestial mechanics.  

The models \eqref{eq: eq motion perturbed Kepler varepsilon=1} and \eqref{eq: eq motion perturbed Kepler} are simple general models which contain the well-studied models of the circular and elliptic restricted three-body problems in a fixed reference frame as special cases, see \cite{OrtegaZhao} and Appendix \ref{Appendix C}. They already present many essential features of the planetary or lunar three-body problems seen as perturbations of two Keplerian elliptic motions. 
Thm. \ref{thm: main} therefore asserts that 
\begin{cor}\label{Cor: main} {In any circular or elliptic restricted three-body problems with masses of the primaries $m_{1}, m_{2}>0$, there exist infinitely many generalized periodic solutions accumulating each of the primaries.} 
\end{cor}
Moreover, Thm. \ref{thm: main} gives an indication that similar phenomenon may as well happen in the non-restricted N-body problems in lunar regimes in which two particles are close to each other while other particles stay sufficiently faraway.  We shall leave this for further investigations. 

We structure this article in the following way: In Section \ref{Sec: 2} we present two regularizations of the problem in the extended phase space, namely extensions of the Levi-Civita regularization for the 2-dimensional case and Moser regularization for all dimensions. We construct proper symplectic coordinates for the regularized systems and discuss the non-degeneracy condition of periodic manifolds of properly-degenerate systems and apply it in our setting to obtain the non-degeneracy in all dimensions in a unified way. In Section \ref{Sec: 3} we present a localization argument which allows us to further localize our analysis near the periodic manifolds $\Lambda_{n}$. In Section \ref{Sec: 4} we further present a rescaling argument which allows us to transform the problem of perturbing infinitely many periodic manifolds into a problem of perturbing a fixed periodic manifold with perturbations with successively decreasing $C^{1}$-norms, which allows us to conclude with Thm. \ref{thm: A} from the Appendix \ref{Appendix A}, established via a local Rabinowitz Floer Homology argument. In Appendix \ref{Appendix B}, we illustrate the construction of symplectic coordinates by proposing variants of the planar Delaunay variables in which the circular motions are regular. In Appendix \ref{Appendix C} we realize circular and elliptic restricted three-body problems as time-periodically forced Kepler problems.



\section{Regularizations of the Forced Kepler Problem in Extended Phase Space}\label{Sec: 2}
We deal with the somehow more general system \eqref{eq: eq motion perturbed Kepler}. By a standard trick we first transform the problem into an autonomous one by  embedding the system into the zero-energy hypersurface of the augmented Hamiltonian $F_{\varepsilon}(p, q, t)+ \tau$ defined on the extended phase space which is the exact symplectic manifold
$$\Bigl(T^{*} (\R^{d} \setminus O)) \times T^{*} S^{1}, d (\sum_{i=1}^{d} p_{i} d q_{i} + \tau d t)\Bigr),$$
 in which $\tau$ denotes the variable conjugate to $t$. The variable $t$ is now considered as a space variable, which for now is identical to the (fictitious) time variable up to a shift by integers.  


We now continue with two different approaches, the first one works for the 2-dimensional case and the second works for all dimensions.

\subsection{Dimension 2: Levi-Civita Regularization}
\subsubsection{Levi-Civita Regularization} We change the fictitious time  $s \mapsto t$ according to $ds/dt=|q|^{-1}$ on the zero-energy level 
$$\{F_{\varepsilon}(p, q, t) + \tau=0\}.$$ 
The resulting system has the same flow as the flow on the zero-energy level of the slowed-down Hamiltonian
$$\{|q|F_{\varepsilon}(p, q, t) + \tau |q|=0\}.$$
We shall use the notation $'=\dfrac{d}{d s}$ to denote the derivative with respect to the new fictitious time $s$. 

In dimension 2 (which contains the 1-dimensional system as a subcase), a regularization is given by further pulling-back the system by the contragradient of the complex square mapping \cite{Goursat}, \cite{Levi-Civita 1904}, \cite{Levi-Civita}:
$$\C \setminus \{0\} \times \C \to \C \times \C, \quad (z, w) \mapsto \Bigl(z^{2}, \dfrac{w}{2 \bar{z}}\Bigr).$$  

The unperturbed regularized Hamiltonian with $\varepsilon=0$ is found to be
$$H_{0}(z, w, \tau, t):=\dfrac{|w|^{2}}{8} + \tau |z|^{2} -1.$$

The regularized Hamiltonian with $\varepsilon \in [0, 1]$ is found to be
$$H_{\varepsilon}(z, w, \tau, t):=\dfrac{|w|^{2}}{8} + \tau |z|^{2} +|z|^{2} \varepsilon U(z^{2}, t, \varepsilon) -1. $$
All these systems now extend regularly near the subset consisting of collisions $\{z=0\}$ of their zero-energy hypersurfaces respectively in the regularized extended phase space $\Bigl(T^{*} \C \times T^{*} S^{1}, d (\Re\{\bar{w} d z\}+ \tau d t)\Bigr)$. The extended systems are called Levi-Civita regularized systems. In the unperturbed regularized system $H_{0}$, the variable $\tau$ is a first integral. When $\tau$ is fixed, $H_{0}$ describes a pair of isotropic harmonic oscillators, \emph{i.e.} a pair of harmonic oscillators with identical frequencies in the $(z, w)$-variables.

\subsubsection{Action-Angle Variables}
We now compute a set of action-angle variables for the unperturbed system given by the Hamiltonian $H_{0}$. Write $z=z_{1}+i z_{2}, w=w_{1}+i w_{2}$ and the Hamiltonian is written as 
$$H_{0}:=(w_{1}^{2}+w_{2}^{2})/8+(z_{1}^{2}+z_{2}^{2})-1.$$
Based on the one-dimensional computation from \cite{Zhao}, we may set
$$z_{1}=2^{-1/4} I_{1}^{1/2} \tau^{-1/4} \cos \theta_{1}, \qquad w_{1}=-2 \cdot 2^{1/4} I_{1}^{1/2} \tau^{1/4} \sin \theta_{1},$$
$$z_{2}=2^{-1/4} I_{2}^{1/2} \tau^{-1/4} \cos \theta_{2},\qquad  w_{2}=-2 \cdot 2^{1/4} I_{2}^{1/2} \tau^{1/4} \sin \theta_{2},$$
\begin{equation}\label{eq: time change}
\tau=\tau, \qquad t=\tilde{t}+4^{-1} \tau^{-1} (I_{1} \sin 2 \theta_{1}+I_{2} \sin 2 \theta_{2}).
\end{equation}
Consequently, we set
$$I_{1}=\dfrac{\mathcal{L}+\mathcal{J}}{2}, I_{2}=\dfrac{\mathcal{L}-\mathcal{J}}{2}, \theta_{1}=\delta+\gamma, \theta_{2}=\delta-\gamma. $$


It is seen that $(\mathcal{L},  \delta, \mathcal{J}, \gamma, \tau, \tilde{t})$ form a set of action-angle variables of $H_{0}$. The pairs $(\mathcal{L}, \delta, \tau, \tilde{t})$ are the fast variables and $H_{0}$ depends only on the fast actions $(\mathcal{L}, \tau)$:
$$H_{0}=\dfrac{\sqrt{2}}{2} \mathcal{L} \sqrt{\tau} -1. $$

We observe that in a neighborhood of $\{H_{0}=0\}$, the pair of fast action variables $(\mathcal{L}, \tau)$ are globally defined, which generates a free Hamiltonian $\T^{2}$-action which we call the fast $\T^{2}$-action. The fast angles $(\delta, \tilde{t})$ are angles associated to this $\T^{2}$-action by our construction of $\mathcal{J}$ and $\gamma$ as long as these latter variables are well-defined. 

The set of variables $(\mathcal{L},  \delta, \mathcal{J}, \gamma, \tau, \tilde{t})$ covers a neighborhood of the zero energy hypersurface $\{H_{0}=0\}$, and in particular each periodic manifold, only open-densely. Indeed it is seen that the variables $(\mathcal{J}, \gamma)$ are action-angle coordinates defined on a dense-open subset of the \emph{orbit space} with fixed $\mathcal{L}$:
$$\Omega_{\mathcal{O}}^{\mathcal{L}} \cong S^{2}.$$
A homeomorphism with $S^{2}$ of the space  $\Omega_{\mathcal{O}}^{\mathcal{L}}$ of planar (possibly circular or rectilinear) ellipses with fixed semi-major axis and with a fixed focus has been explained e.g. in \cite{Albouy}. 

We now argue that we may cover a neighborhood of our energy hypersurface with a finite set of similar coordinate charts. 

Indeed, proceeding with a symplectic reduction with respect to the fast Hamiltonian $\T^{2}$-action, we see that the symmetric group $SO(3)$ of the system of two isotropic harmonic oscillators acts in a Hamiltonian way on the orbit space  $\Omega_{\mathcal{O}}^{\mathcal{L}}$, which in this case corresponds to the natural action of $SO(3)$ on $S^{2}$. This then restricts to the Hamiltonian action of any $S^{1}$-subgroup of $SO(3)$. This last $S^{1}$-action is everywhere free on $\Omega_{\mathcal{O}}^{\mathcal{L}}$ outside two antipodal points, and extends the fast $\T^{2}$-action to a Hamiltonian $\T^{3}$-action on the regularized extended phase space. We may then apply Arnold-Liouville theorem to the Lagrangian torus fibration associated to this Hamiltonian $\T^{3}$-action which then gives us a set of action-angle variables locally around any Lagrangian torus within this fibration. By passing from the action-angle variables locally to Darboux coordinates for variables on $\Omega_{\mathcal{O}}^{\mathcal{L}}$, we get a set of well-defined variables $(\mathcal{L},  \delta, \tau, \tilde{t}, \xi, \zeta)$, where $(\xi, \zeta)$ are Darboux coordinates defined in an open set $\Omega_{\mathcal{O}}^{loc, \mathcal{L}}$ of $\Omega_{\mathcal{O}}^{\mathcal{L}}$ which may be assumed to be small. Note that the angles $(\delta, \tilde{t})$ may have different zero-sections when compared to their analogues in the set of variables $(\mathcal{L},  \delta, \mathcal{J}, \gamma, \tau, \tilde{t})$ that we have previously constructed on overlap of charts, but this will not be a concern to us and we shall keep their notations unchanged in different charts.
By Taking images of $\Omega_{\mathcal{O}}^{loc, \mathcal{L}}$ under the $SO(3)$-action we see that, since $\Omega_{\mathcal{O}}^{\mathcal{L}}$ is compact on which $SO(3)$ acts transitively, there exists a finite cover of $\Omega_{\mathcal{O}}^{\mathcal{L}}$ by the $SO(3)$-images of $\Omega_{\mathcal{O}}^{loc, \mathcal{L}}$, which then gives rise to finitely many charts of a neighborhood of the zero-energy hypersurface of $\{H_{0}=0\}$ in the extended phase space by the $SO(3)$-invariance of the system. 
 

We now aim at  better understanding the dependence of $\xi, \zeta$ on $\mathcal{L}$. By the $SO(3)$-invariance it is enough to understand this in a particular chart. But we see that with the initial set of variables $(\mathcal{L},  \delta, \mathcal{J}, \gamma, \tau, \tilde{t})$, the quantities $(\mathcal{J}/\mathcal{L}, \gamma)$ are functions defined on $\Omega_{\mathcal{O}}=\Omega_{\mathcal{O}}^{1}$ (the space $\Omega_{\mathcal{O}}^{\mathcal{L}}$ with $\mathcal{L}=1$) independent of $\mathcal{L}$ and $\tau$. Passing to Darboux coordinates and with the $SO(3)$-invariance we conclude that any set of Darboux coordinates $(\xi, \zeta)$ on the normalized orbit space $\Omega_{\mathcal{O}}=\Omega_{\mathcal{O}}^{1}$ gives rise to a set $(\mathcal{L},  \delta, \mathcal{J}, \gamma, \tau, \tilde{t}, \hat{\xi}, \hat{\zeta})$ for $(\hat{\xi}, \hat{\zeta})=(\sqrt{\mathcal{L}} \xi, \sqrt{\mathcal{L}} \zeta)$. 

We summarize this analysis in the following Proposition:


\begin{prop} Any Darboux coordinates $(\xi, \zeta)$ of the space  $(\Omega_{\mathcal{O}}, \omega_{0})$ gives rise to a set of local symplectic coordinates $(\mathcal{L},  \delta, \tau, \tilde{t}, \hat{\xi}, \hat{\zeta})$ of $T^{*} \C \times T^{*} S^{1}$ where $(\hat{\xi}, \hat{\zeta})=(\sqrt{\mathcal{L}} \xi, \sqrt{\mathcal{L}} \zeta)$. A neighborhood of the energy hypersurface $\{H_{0}=0\}$ in the extended phases space is covered by finitely many such symplectic charts.
\end{prop}

\subsubsection{Periodic Manifolds of the Regularized Kepler Problem in Extended Phase Space}

For a Hamiltonian system defined on an exact symplectic manifold 
$$(M, \omega=d \lambda, H)$$
 we denote by $\phi_{t}$
its flow. We let the set 
$\hbox{Per}_{0}^{H} \subset M \times \R$  be the set of points $(m, \eta) \in M \times \R$ such that 
$$H(m)=0, d H(m) \neq 0, \phi_{\eta}(m)=m. $$
It follows that for $(m, \eta) \in \hbox{Per}_{0}^{H}$, the orbit $\phi_{t} (m)$ is a periodic orbit of the system in $H^{-1}(0)$ with period $\eta$.

Following \cite{Weinstein}, a subset $\Lambda \subset  \hbox{Per}_{E}^{H}$ is called a periodic manifold if $\Lambda$ is a submanifold of $M \times \R$ and if the restriction on $\Lambda$ of the projection $\pi: M \times \R \mapsto M$ is an embedding. A periodic manifold is called \emph{non-degenerate} if at any $(m, \eta) \in \Lambda$ there holds
\begin{equation}\label{Condition: Nondegeneracy}
T_{m} \pi(\Lambda)=(D \phi_{\eta}-Id)^{-1} (\hbox{span}\{X_{H}(m)\}). 
\end{equation}
It is under this non-degeneracy hypothesis that Weinstein's theorem is applicable.

It is direct to check, see e.g.  \cite{Boscaggin--Ortega--Zhao}, that there always holds 
$$T_{m} \pi(\Lambda) \subset (D \phi_{\eta}-Id)^{-1} (\hbox{span}\{X_{H}(m)\}).$$ 
Therefore the non-degeneracy condition is equivalent to state that the two spaces $T_{m} \pi(\Lambda)$ and $(D \phi_{\eta}-Id)^{-1} (\hbox{span}\{X_{H}(m)\})$ have the same dimension. 


We denote by $\Lambda_{n}$ the periodic manifolds of the regularized Kepler problem $H_{0}$ in extended phase space with prime period $1/n$ with respect to the initial time variable $t$. 

Observe that for a periodic orbit, a period for $t$ is automatically a period for $\tilde{t}$ as well. We may therefore assume that $\tilde{t}$ has prime period $1/n$, with correspondingly $\mathcal{L}=\mathcal{L}_{n}, \tau=\tau_{n}$. A direct computation shows that the minimal period for the angle $\delta$ to go through $ \pi$ (instead of $2 \pi$, since the Levi-Civita regularization mapping is two-to-one) is $S_{n}= \pi \sqrt{2}\tau_{n}^{-1/2}$. Therefore by assumption
$$S_{n} \cdot \tilde{t}'(\mathcal{L}_{n}, \mathcal{\tau}_{n})= \pi \sqrt{2}\tau_{n}^{-1/2}\cdot \dfrac{\sqrt{2}}{4} \mathcal{L}_{n} /\sqrt{\tau_{n}}=\dfrac{1}{2} \pi \mathcal{L}_{n} \tau_{n}^{-1}=1/n.$$ 
Together with the energy condition 
$$\dfrac{\sqrt{2}}{2} \mathcal{L}_{n} \sqrt{\tau}_{n}-1=0$$
 we obtain
\begin{equation}\label{eq: Ln, taun}
\mathcal{L}_{n}=2^{2/3} \pi^{-1/3}  n^{-1/3}, \qquad \tau_{n}=2^{-1/3} \pi^{2/3}  n^{2/3},
\end{equation}
and
$$S_{n}=(2 \pi)^{2/3} n^{-1/3}.$$

In the planar case, these manifolds have been analyzed in \cite{Boscaggin--Ortega--Zhao} in which it is verified that they are non-degenerate. 

\subsubsection{Non-Degeneracy of Integrable Hamiltonian in Action-Angle Form in Extended Phase Space } For the purpose of uniformly treating non-degeneracy conditions we now revisit this non-degeneracy with the help of the action-angle forms of an integrable (possibly properly-degenerate) Hamiltonian. 
We consider an N-degrees of freedom Hamiltonian function in extended phase space of the form
$$H(I_{1}, \cdots, I_{n}, \tau)=0,$$
in which $n \le N$, within the chart given by the symplectic variables 
$$(I_{1}, \theta_{1}, \cdots, I_{n}, \theta_{n}, \tau, t, \hat{\xi}_{n+1}, \hat{\zeta}_{n+1}, \cdots, \hat{\xi}_{N}, \hat{\zeta}_{N})$$
containing the basepoint point $m$, such that the first $n+1$-pairs are (partial) action-angle variables. The fact that the Hamiltonian is independent of all the variables $(\hat{\xi}_{1}, \hat{\zeta}_{1}, \cdots, \hat{\xi}_{n}, \hat{\zeta}_{n})$ implies that the entire periodic orbits through $m$ is contained in this symplectic chart. In particular, it is enough to investigate the non-degeneracy condition in each of these charts.

For a solution to be periodic, the frequency vector 
$$(\nu_{1}, \cdots, \nu_{n}, \nu_{\tau}):=\Bigl(\dfrac{\partial H(I_{1}, \cdots, I_{n}, \tau)}{\partial I_{1}}, \cdots, \dfrac{\partial H(I_{1}, \cdots, I_{n}, \tau)}{\partial I_{n}}, \dfrac{\partial H(I_{1}, \cdots, I_{n}, \tau)}{\partial \tau}\Bigr)$$
need to satisfy that  for all $i=1, \cdots, n$, the corresponding frequency $\nu_{i}$ is an integer multiple of $\nu_{\tau}$ and  $\nu_{\tau}\neq 0$. This latter implies in particular that we may at least locally express $\tau=\tau(I_{1}, \cdots, I_{n})$ as a function of $(I_{1}, \cdots, I_{n})$ by implicit function theorem. 

We now consider a periodic manifold of the system given by the constraints $(I_{1}=I_{1}^{0}, \cdots, I_{n}=I_{n}^{0},\tau=\tau_{0})$. 
We set 
$$\nu_{\tau}^{0}=\dfrac{\partial H(I_{1}, \cdots, I_{n}, \tau)}{\partial \tau} (I_{1}=I_{1}^{0}, \cdots, I_{n}=I_{n}^{0}).$$

The tangent space at a point is of dimension $2 N -n +2$ and is given by 
$$\hbox{span}\{\partial_{\theta_{1}}, \cdots \partial_{\theta_{n}}, \partial_{t}, \partial_{\hat{\xi}_{n+1}}, \partial_{\hat{\zeta}_{n+1}}\cdots  \partial_{\hat{\xi}_{N}}, \partial_{\hat{\zeta}_{N}}\}.$$

The flow direction is given by the vector
$$\dfrac{\partial H}{\partial {I_{1}}}(I_{1}^{0}, \cdots, I_{n}^{0}, \tau^{0}) \partial {\theta_{1}} + \cdots + \dfrac{\partial H}{\partial {I_{n}}}(I_{1}^{0}, \cdots, I_{n}^{0}, \tau^{0}) \partial {\theta_{n}}+\nu_{\tau}^{0}\partial_{t}.$$

Moreover we have $\eta=(\nu_{\tau}^{0})^{-1}$ and the mapping $\phi_{\eta}$ takes the form
\begin{eqnarray*}
\phi_{\eta}: &(I_{1}, \cdots I_{n}, \tau,, \theta_{1}, \cdots \theta_{n}, t, \hat{\xi}_{n+1}, \hat{\zeta}_{n+1} \cdots \hat{\xi}_{N}, \hat{\zeta}_{N}) \\
                 \mapsto & (I_{1}, \cdots I_{n}, \tau, \theta_{1}+ ({\nu_{\tau}^{0}})^{-1} \dfrac{\partial H}{\partial {I_{1}}}(I_{1}, \cdots, I_{n}, \tau) , \cdots \theta_{n} + ({\nu_{\tau}^{0}})^{-1} \dfrac{\partial H}{\partial {I_{n}}}(I_{1}, \cdots, I_{n}, \tau),\\ &t+({\nu_{\tau}^{0}})^{-1} \dfrac{\partial H}{\partial \tau}(I_{1}, \cdots, I_{n}, \tau), \hat{\xi}_{n+1}, \hat{\zeta}_{n+1}\cdots \hat{\xi}_{N}, \hat{\zeta}_{N}).
\end{eqnarray*}

Thus in matrix form we may write
\begin{eqnarray*}
D \phi_{\eta}: &\Bigl(Id_{(2n+2) \times (2n+2)}  +  \begin{pmatrix} 0 & 0  \\ ({\nu_{\tau}^{0}})^{-1} \hbox{Hess}_{(I_{1}, \cdots I_{n}, \tau)} H(I_{1}^{0}, \cdots, I_{n}^{0}, \tau_{0}) & 0\end{pmatrix} \Bigr) \\ & \times Id_{(2 N-2 n) \times (2 N-2 n)}.
\end{eqnarray*}
where we have denoted by
$$\hbox{Hess}_{(I_{1}, \cdots I_{n}, \tau)} H (I_{1}^{0}, \cdots, I_{n}^{0},\tau_{0})$$
the Hessian matrix of $H$ with respect to $(I_{1}, \cdots, I_{n}, \tau)$ evaluated at the periodic manifold $(I_{1}=I_{1}^{0}, \cdots, I_{n}=I_{n}^{0}, \tau=\tau_{0})$.
In view of \eqref{Condition: Nondegeneracy} we may now conclude:

\begin{prop}\label{prop: non-degeneracy} The periodic manifold $(I_{1}=I_{1}^{0}, \cdots, I_{n}=I_{n}^{0})$ is non-degenerate if and only if the matrix
$$ \hbox{Hess}_{I_{1}, \cdots I_{n}, \tau} H (I_{1}^{0}, \cdots, I_{n}^{0},\tau_{0})$$
is non-degenerate. 
\end{prop}
We remark that in the non-properly-degenerate case, this non-degeneracy condition corresponds to the Kolmogorov non-degeneracy condition well-known in KAM theory. 

Applying this Proposition to the Levi-Civita regularized Kepler Hamiltonian
$$H_{0}=\dfrac{\sqrt{2}}{2} \mathcal{L} \sqrt{\tau} -1, $$
we directly obtain the fact that 

\begin{prop} The periodic manifolds $\Lambda_{n}$ in the $d$-dimensional Moser regularized Kepler problem are non-degenerate for all $d \in \Z_{+}$.
\end{prop}



\subsubsection{Computation of Action Values of Periodic Manifolds}\label{Sec: Action Values}






On $\{H_{0}=0\}$, the action of a corresponding periodic solution with prime period $1/n$ for a period with length $1$ (both refer to the initial $t$-variable) is computed via the integral 
\begin{eqnarray*}
&n \int_{0}^{S_{n}} (\mathcal{L}_{n} \delta'(\mathcal{L}_{n},\tau_{n})+\tau_{n} \tilde{t}'(\mathcal{L}_{n},\tau_{n})) ds \\
&=n \int_{0}^{2 \pi \sqrt{2}\tau_{n}^{-1/2}} (\mathcal{L}_{n} \dfrac{\sqrt{2}}{2} \sqrt{\tau}_{n}+\tau_{n} \dfrac{\sqrt{2}}{4} \mathcal{L}_{n} /\sqrt{\tau}_{n}) ds.
\end{eqnarray*}
  
Now substitute the corresponding values of $\mathcal{L}_{n}$ and $\tau_{n}$ from \eqref{eq: Ln, taun}, we find directly that
$$\mathcal{A}_{0}^{(n)}=3 \cdot 2^{-1/3} \pi^{2/3} n^{2/3}=\mathcal{A}_{0}^{(1)} n^{2/3}.$$

From this we may state the following observations:
\begin{prop}\label{prop: action value}
\begin{itemize}
It holds that 
\item $\mathcal{A}_{0}^{(n)}$ is monotone in $n$;
\item $\lim_{n \to \infty} \mathcal{A}_{0}^{(n)}=\infty$;
\item the gaps between the action values of $\Lambda_{n+1}$ and $\Lambda_{n}$ is
$$\mathcal{A}_{0}^{(1)} ((n+1)^{2/3}-n^{2/3}) 
\sim n^{-1/3}.$$
\end{itemize}
\end{prop}



\subsection{{All Dimensions: Moser Regularization}}
\subsubsection{Moser Regularization}We recall a regularization of the Kepler problem and its regular perturbations in all dimensions now typically attributed to Moser \cite{Moser}. 

We start with the Kepler problem 
$$\Bigl(T^{*} (\R^{d} \setminus \{0\}), \omega_{0}, \dfrac{\|p\|^{2}}{2}-\dfrac{1}{\|q\|}\Bigr)$$
with coordinates $(q, p) \in \R^{d} \setminus \{0\} \times \R^{d}$. On the energy level $-1/2$ we may change time and arrive at the  zero-energy level of the system
$$\Bigl(T^{*} (\R^{d} \setminus \{0\}), \omega_{0}, \dfrac{\|q\|(\|p\|^{2}+1)}{2}-1=0\Bigr).$$
Proceeding with the canonical change of variables 
$$\R^{d} \times \R^{d} \to \R^{d} \mapsto \R^{d}, (-x, y) \mapsto (p,q),$$ 
we arrive at the following system on its energy level $1$
$$\Bigl\{\dfrac{\|y\|(\|x\|^{2}+1)}{2}=1\Bigr\}.$$ 
Following Moser we may treat this as the stereographic projection of the system
$$(T^{*}_{+} (S^{d} \setminus N), \omega_{0}, \|v\|=1),$$
where $T^{*}_{+} S^{d}$ denotes the positive cotangent bundle in $T^{*} (S^{d})$, and $v$ denotes the momenta variables. We see that this is equivalent to the system of geodesic flow
$$\Bigl(T^{*}_{+} (S^{d} \setminus N), \omega_{0}, \dfrac{\|v\|^{2}}{2}=\dfrac{1}{2} \Bigr).$$
The set of collisions, or equivalently, the set corresponding to infinite velocities, are now transformed to the fibre $T^{*}_{N} S^{d}$. Moser's regularization is completed by the further addition of this fibre and the system thus defined
$$\Bigl(T^{*}_{+} S^{d}, \omega_{0}, \dfrac{\|v\|^{2}}{2}=\dfrac{1}{2} \Bigr)$$
is regular with great circle orbits with energy $1/2$. A particle moves on such a great circle with unit velocity, therefore has period $2 \pi$. 

By properly rescaling the constant we arrive at a regularization for all negative energies. We now work out more details. 

We take a $d$-dimensional centered sphere $S_{r}^{d}$ with radius $r$ in $\R^{d+1}$ and consider the plane $\R^{d} \times \{0\}$ to be the plane to be projected to from the north pole $\mathcal{N}=(0,\cdots,0,r)$ of the sphere. Let $$u=(u_{1}, \cdots, u_{d+1}) \in S_{r}^{d} \setminus \mathcal{N} $$ be a point on this sphere different from $\mathcal{N}$. We let $v=(v_{1}, \cdots, v_{d+1}) \in \R^{d+1}$ be a (co-)tangent vector at this point, so explicitly we have
$$u \neq \mathcal{N}, u\cdot u=r^{2}, u \cdot v=0.$$
We now let $x=(x_{1}, \cdots, x_{d}) \in \R^{d}$ be such that $u$ projects to $(x, 0)$ by stereographic projection from $\mathcal{N}$ and $y=(y_{1}, \cdots, y_{d})$ be a (co-)vector at the point $x$, such that the 1-form 
$$\sum_{j=1}^{d+1} v_{j} d u_{j} =\sum_{i=1}^{d} y_{i} d x_{i}$$ 
is preserved.

We therefore have
\begin{equation}\label{eq: stereographic 1}
x_{i}=\dfrac{r u_{i}}{r-u_{d+1}}, y_{i}=\dfrac{r-u_{d+1}}{r} v_{i} + \dfrac{u_{i} v_{d+1}}{r}, i=1,\cdots, d, 
\end{equation}
with inverse
\begin{eqnarray*}
& u_{i}=\dfrac{2 r^{2} x_{i}}{r^{2}+\|x\|^{2}}, u_{d+1}=\dfrac{\|x\|^{2}-r^{2}}{\|x\|^{2}+r^{2}}r, \\ & v_{i}=\dfrac{\|x\|^{2}+r^{2}}{2 r^{2}} y_{i} - \dfrac{(x \cdot y) x_{i}}{r^{2}}, v_{d+1}=\dfrac{x \cdot y}{r}, i=1,\cdots, d.
\end{eqnarray*}
From these equalities it follows that
$$\|v\|=\Bigl(\dfrac{\|x\|^{2}+r^{2}}{2 r^{2}}\Bigr) \|y\|.$$

Therefore for the Kepler problem, once we fix a negative energy $-f<0$, after a proper time reparametrization, we arrive at the system
$$\Bigl(T^{*}_{+} (\R^{d} \setminus \{0\}), \omega_{0}, 2 f\dfrac{\|q\|(\|p\|^{2}+2 f)}{2 \cdot(2f)}-1=0\Bigr),$$
which after the symplectic change of variables $(-x, y) \mapsto (p,q)$ is the stereographic projection of the system
$$(T^{*}_{+} (S^{d}_{\sqrt{2 f}} \setminus \{\mathcal{N}\}), \omega_{0}, 2 f \|v\|=1).$$
This is equivalent to the reparametrized geodesic flow system
$$\Bigl(T^{*} (S^{d}_{\sqrt{2 f}} \setminus \{\mathcal{N}\}), \omega_{0}, \dfrac{(2 f \|v\|)^{2}}{2}=\dfrac{1}{2}\Bigr),$$
so the orbits are great circles with radius $\sqrt{2 f}$ and the velocities have norm $1/2f$. The period is therefore $2 \pi (\sqrt{2 f})^{3}=4 \sqrt{2} \pi f^{3/2}$. 

Unlike Levi-Civita and  Kustaanheimo-Stiefel regularizations as have been used in \cite{Boscaggin--Ortega--Zhao}, \cite{Verzini-Ortega-Barutello}, this regularization does not create additional symmetry of the regularized system and is valid for all dimensions in a uniform way. On the other hand, here we have to work on curved space instead of Euclidean spaces.
In \cite{Kummer}, the relationship between the Levi-Civita/Kustaanheimo-Stiefel regularization and the Moser regularization has been established, and Moser regularization can be regarded in the planar case as a quotient of Levi-Civita regularization by a $\Z_{2}$-symmetry, and in the spatial case as a symplectic quotient of the Kustaanheimo-Stiefel regularization by a Hamiltonian $S^{1}$-symmetry respectively. 

\subsubsection{Extension to Forced Kepler Problem}\label{Subsubsection: Extension to Forced Kepler problem}
Moser regularization extends to forced Kepler problem in the extended phase space in the following way: Starting from the system
$$\Bigl(T^{*} (\R^{d} \setminus \{0\}) \times T^{*} S^{1}, \omega_{0}, \dfrac{\|p\|^{2}}{2}-\dfrac{1}{\|q\|} + \varepsilon U(q, t, \varepsilon)+\tau=0\Bigr)$$
with coordinates $(q, p, t, \tau) \in \R^{d} \setminus \{0\} \times \R^{d}  \times S^{1} \times \R$, we change time
 so that
we obtain the system
$$\Bigl(T^{*} (\R^{d} \setminus \{0\}), \omega_{0}, \dfrac{\|q\|(\|p\|^{2}+2\tau)}{2}+ \varepsilon \|q\| U(q, t, \varepsilon)-1=0\Bigr).$$
with energy-level zero. We set $(x, y) = (-p,q)$ as in previous discussions.
 
It follows from \eqref{eq: stereographic 1} that
$$\|p\|^{2}=\|x\|^{2}=\dfrac{2 \tau (2 \tau-u_{d+1}^{2})}{(\sqrt{2 \tau}-u_{d+1})^{2}},$$
and therefore
$$\dfrac{2\tau}{\|p\|^{2}+2\tau}=\dfrac{\sqrt{2 \tau}-u_{d+1}}{2 \sqrt{2 \tau}},$$
which extends to a  regular function on the sphere $S^{d}_{\sqrt{2 \tau}}$. Therefore 
\begin{equation}\label{eq: q from u, v}
\|y\|=\|q\|=\|v\| \dfrac{\sqrt{2 \tau}-u_{d+1}}{\sqrt{2 \tau}}. 
\end{equation}

Therefore by stereographic projection the system extends to a system
$$(\{T^{*}_{+} (S^{d}_{\sqrt{2 \tau}}):(\tau, t) \in T^{*} S_{1}\}, \omega_{0}, 2 \tau \|v\|-1+ \varepsilon \|v\|\dfrac{\sqrt{2 \tau}-u_{d+1}}{ \sqrt{2 \tau}}  U(q(v, u), t, \varepsilon) =0)$$
for $u \in S^{d}_{\sqrt{2 \tau}}$. The space is a (trivial) fibre bundle with fibres $T^{*}_{+} (S^{d}_{\sqrt{2 \tau}})$ over $T^{*} S^{1}$, with fibres depending on $\tau$ and independent of $t$. 


To normalize the fibres to be $\tau$-independent we make a symplectic rescaling
$$v=\tilde{v} /\sqrt{2 \tau} , \qquad u=\sqrt{2 \tau} \tilde{u},$$
and accordingly take a shift of the $t$-variable
$$\tilde{t}=t-\phi(\tau, \tilde{u}, \tilde{v})$$
for some function $\phi$ for the purpose of keeping the change of variables canonical. The function $\phi$ is implicitly determined by the relationship {(after lifting $t, \tilde{t}$ to be defined on $\R$)}:
\begin{equation}\label{eq: t change}
{\sum_{i} v_{i} d u_{i}-t d \tau=\sum_{i} \tilde{v}_{i} d \tilde{u}_{i}-\tilde{t} d \tau.}
\end{equation}
We shall only need the form of $\phi$ in local coordinates which we shall determine later. 

We thus have 
\begin{equation}\label{exp: q}
\|q\|=\|\tilde{v}\| \tau^{-1/2} (1-u_{d+1}).
\end{equation}

The resulting system is equivalent to
$$(T^{*}_{+} S^{d}_{1}  \times T^{*} S^{1}, \omega_{0}, \sqrt{2 \tau} \|\tilde{v}\|-1+\varepsilon \|\tilde{v}\| (1-\tilde{u}_{d+1}) (2 \tau)^{-1/2} U(q(\tilde{v}/\sqrt{2 \tau}, u),\tilde{t}, \varepsilon) =0).$$

{When $\varepsilon=0$, we see that the corresponding energy hypersurface of the regularized Kepler problem is $T^{*}_{(2 \tau)^{-1/2}} S^{d}_{1}  \times T^{*} S^{1}$. }

{In the above system the perturbation term takes the form
$$\varepsilon U=\varepsilon \|\tilde{v}\| (1-\tilde{u}_{d+1}) (2 \tau)^{-1/2} U(q(\tilde{v}/\sqrt{2 \tau}, u),\tilde{t}, \varepsilon).$$ 
In view of \eqref{eq: q from u, v}, \eqref{eq: estimate on U} and the energy constraint,  by passing to a local chart of the sphere $S^{d}$ around the north pole given by $(\tilde{u}_{1}, \cdots, \tilde{u}_{d})$ and using the identity
$$1-\tilde{u}_{d+1}=\dfrac{\tilde{u}_{1}^{2}+\cdots+\tilde{u}_{d}^{2}}{1+\tilde{u}_{d+1}}$$
 we get that in a sufficiently small neighborhood $\{\tilde{u}_{1}^{2}+\cdots+\tilde{u}_{d}^{2} \le \tilde{\varepsilon}\},$ for some sufficiently small, $0<\tilde{\varepsilon}<<1$ of the north pole, and for $\tau \ge \tau_{*}>0$, the perturbation term $\varepsilon U$ has the estimate
 $$\|\varepsilon U\| \le \tilde{C}_{0}  (\tilde{u}_{1}^{2}+\cdots+\tilde{u}_{d}^{2})^{2},$$
for some constant $\tilde{C}_{0}$ which depending only on $\tau_{*}$ and $\tilde{\varepsilon}$. 

From these discussion we conclude
\begin{prop}\label{prop: patching}
The perturbation term $\varepsilon U$ is regular near the set of collisions and vanishes on this set together with its first, second and third derivatives. 
\end{prop}}

\subsubsection{A Set of Symplectic Coordinates}\label{Sec: 2.2.3}
As in the two-dimensional case, the space $T^{*}_{+} S^{d}_{1}  \times T^{*} S^{1}$ carries a free Hamiltonian $\T^{2}$-action associated to the pair of variables $(\mathcal{I}, \tau)$, in which $\mathcal{I}$ is a moment map of the Hamiltonian $S^{1}$-action on $T^{*}_{+} S^{d}_{1}$ given by the geodesic flow on $S^{d}_{1}$. 

When we restrict the system $H_{0}$ with any fixed $\tau>0$ to the cotangent bundle $T^{*} S^{1}$ of a great circle $S^{1}$ in $S^{d}$, we observe that both $\mathcal{I}$ and $\|\tilde{v}\|$ generates the same $S^{1}$-action on $T^{*} S^{1}$, thus in this case $\mathcal{I}$ and $\|\tilde{v}\|$ agree up to an additive constant, which may be discarded by shifting $\mathcal{I}$ by this constant. We may thus set $\mathcal{I}=\|\tilde{v}\|$ by the $SO(d)$-invariance of the system $H_{0}$ . Therefore we have 
$$H_{0}=\sqrt{2 \tau} \mathcal{I} -1.$$

The symplectic quotient of the space $T^{*}_{+} S^{d}_{1}  \times T^{*} S^{1}$ by this free Hamiltonian $\T^{2}$-action at a level with fixed $\mathcal{I}$ and with total energy zero is the space $\Omega_{\mathcal{O}}^{\mathcal{I}} \cong T_{1}^{*} S^{d}/S^{1}$ that we refer to as the orbit space. The action of the symmetric group $SO(d+1) \times S^{1}$ on  $T^{*}_{+} S^{d}_{1}  \times T^{*} S^{1}$ descends to an action of $SO(d+1)$ on the orbit space $\Omega_{\mathcal{O}}^{\mathcal{I}}$ which is Hamiltonian. Choosing any maximal torus $\mathcal{T} \subset SO(d+1)$ we may now apply Arnold-Liouville theorem to deduce the existence of a set of local symplectic coordinates $(\mathcal{I}, \theta, \tau, \tilde{t}, \hat{\xi}_{1}, \hat{\zeta}_{1}, \cdots, \hat{\xi}_{d}, \hat{\zeta}_{d})$. The coordinates $(\hat{\xi}_{1}, \hat{\zeta}_{1}, \cdots, \hat{\xi}_{d}, \hat{\zeta}_{d})$ are defined on an open subset $\Omega_{\mathcal{O}}^{\mathcal{I}, loc} \subset \Omega_{\mathcal{O}}^{\mathcal{I}}$ and by compactness of $\Omega_{\mathcal{O}}^{\mathcal{I}}$ and transitivity of the $SO(d+1)$ we conclude that we may cover $\Omega_{\mathcal{O}}^{\mathcal{I}}$  with finitely many images of $\Omega_{\mathcal{O}}^{\mathcal{I}, loc}$ under the action of $SO(d+1)$.


In the case $d=2$, by comparing both forms of the regularized Hamiltonians by Moser and Levi-Civita regularizations respectively and viewing the 2-to-1 cover of the angles, we find that
$$\mathcal{I}=\mathcal{L}/2, \qquad \theta=2 \delta.$$

We denote by $\Omega_{\mathcal{O}}=\Omega_{\mathcal{O}}^{1}$ and call it the normalized orbit space. As in the $d=2$ case with Levi-Civita regularization, we have the following: 

\begin{prop}Any Darboux coordinates $(\xi_{i}, \zeta_{i})$ on $(\Omega_{\mathcal{O}}, \omega_{0})$ and the symplectic form of  $(T^{*}_{+} (S^{d}_{1} \setminus \{\mathcal{N}\}) \times T^{*} (\R/ \Z)$ can be written  as 
$$d \mathcal{I} \wedge d \theta + d \tau \wedge d \tilde{t} + \sum_{i=1}^{d-1} d (\sqrt{ \mathcal{I}} \xi_{i}) \wedge d (\sqrt{ \mathcal{I}} \zeta_{i})=d \mathcal{I} \wedge d \theta + d \tau \wedge d \tilde{t} +  \sum_{i=1}^{d-1} d \hat{\xi}_{i} \wedge d  \hat{\zeta}_{i},$$
for $(\hat{\xi}_{i}, \hat{\zeta}_{i})=(\sqrt{\mathcal{I}} \xi_{i},\sqrt{\mathcal{I}} \zeta_{i})$. 
\end{prop}


In this set of Darboux coordinates we have from \eqref{exp: q} that
$$q_{i} =\mathcal{I} \tau^{-1/2} f_{i}(\theta, \xi_{i}, \zeta_{i})$$ 
for some regular functions $f_{i}$ independent of $\mathcal{I}$ and $\tau$. Consequently we may write
$$|q|U(q, t)=\mathcal{I}^{2} \tau^{-1} \tilde{U}(\mathcal{I}, \theta, \xi_{i}, \zeta_{i}, t),$$
where by \eqref{eq: estimate on U}, the function $ \tilde{U}(\mathcal{I}, \theta, \xi_{i}, \zeta_{i}, t)$ satisfies
$$ \tilde{U}( \mathcal{I}, \theta, \xi_{i}, \zeta_{i}, t)= \hat{U}( \theta, \xi_{i}, \zeta_{i}, t) + O(\mathcal{I} \tau^{-1/2}).$$

Moreover, since the functions $\tilde{t}, t, \mathcal{I}, \tau$ are $SO(d)$-invariant, by comparing to the formula \eqref{eq: time change} in the 2-dimensional case, we conclude in the same way that the change of time variable takes the form
\begin{equation}\label{eq: time change multidimensional}
\tilde{t}=t+\mathcal{I} \tau^{-1} \tilde{\phi}(\theta, \xi_{i}, \zeta_{i})
\end{equation}
for some smooth function $\tilde{\phi}$.

To be in consistence in the 2-dimensional case with Levi-Civita regularization, we may as well further set 
$$\mathcal{I}=\mathcal{L}/2, \theta=2 \delta.$$

The unperturbed Hamiltonian takes the form 
$$H_{0}=\dfrac{\sqrt{2}}{2} \mathcal{L} \sqrt{\tau} -1$$
which is the same as the unperturbed Levi-Civita regularized system in dimension 2. This expression is now valid in all dimensions.

\subsubsection{Periodic Manifolds of the Moser-regularized Kepler Problem and Their Non-Degeneracy}
In terms of the variables $(\mathcal{L}, \tau)$, the periodic manifold $\Lambda_{n}$ of the Moser-regularized Kepler problem with prime period $1/n$ is characterized by the condition
$$\{\mathcal{L}=\mathcal{L}_{n}, \tau=\tau_{n}\}.$$
Since all the change of variables are canonical and the action value is constant on $\Lambda_{n}$. By restricting to the 2-dimensional case, we conclude that the corresponding action value is $\mathcal{A}_{0}^{(n)}$ which has been computed in Section \ref{Sec: Action Values}.  

Applying Prop. \ref{prop: non-degeneracy} we conclude that
\begin{prop} For any $d \in \Z_{+}$, the periodic manifolds $\Lambda_{n}$ of the Moser-regularized Kepler problem in extended phase space are non-degenerate.
\end{prop}



\section{Localization in Space}\label{Sec: 3}

In this section we shall explain a localization procedure of the regularized system. The argument is presented with the Levi-Civita regularized system in dimension 2 for fixing the ideas, but it works for the Moser regularized system in any dimension without much change, which we shall explain at the end of this section.

\subsection{The L.C. Regularized Systems in Action-Angle Form}

From previous discussions we have seen that the regularized extended Hamiltonian is written in a proper chart as 
\begin{equation}\label{eq: reg Ham}
H_{\varepsilon}(\mathcal{L}, \mathcal{\delta}, \xi, \zeta, \tau, t, \varepsilon):=  \dfrac{\sqrt{2}}{2} \mathcal{L}\sqrt{\tau}-1+ \varepsilon \mathcal{L}^{2} \tau^{-1} \cdot \tilde{U}(\mathcal{L}, \mathcal{\delta}, \xi, \zeta, t, \varepsilon),
\end{equation}
in which the function 
$$ \tilde{U}(\mathcal{L}, \mathcal{\delta}, \xi, \zeta, t, \varepsilon) $$ 
is $C^{\infty}$ regular in all of its variables and depends continuously on the parameter $\varepsilon \in [0, 1]$. Moreover, we have that 
$$ \tilde{U}(\mathcal{L}, \mathcal{\delta}, \xi, \zeta, t, \varepsilon)= \hat{U}(\mathcal{\delta}, \xi, \zeta,\tau, t, \varepsilon) + O(\mathcal{L} \tau^{-1/2}). $$
A priori the terms in $O(\mathcal{L} \tau^{-1/2})$ may depend on all of the listed variables.

In order to work within the corresponding symplectic chart
$(\mathcal{L}, \mathcal{\delta}, \hat{\xi}, \hat{\zeta}, \tau, t),$
we shall treat $\xi, \zeta, t$ as functions defined in this chart given by the following expressions
\begin{equation}\label{eq: variable change LC}
\xi=\mathcal{L}^{-1/2} \hat{\xi},\quad \zeta=\mathcal{L}^{-1/2} \hat{\zeta}, \quad t =\tilde{t}+\mathcal{L} \tau^{-1} \tilde{\phi} (\xi, \zeta),
\end{equation}
in which $\tilde{\phi}$ is a regular function of $(\xi, \zeta)$ as defined in \eqref{eq: time change}.


Recall that we use $'=\dfrac{d}{ds}$ to denote the time derivative of a quantity with respect to the reparametrized fictitious time $s$. 

The non-trivial Hamiltonian equations associated to the unperturbed system $H_{0}=\dfrac{\sqrt{2}}{2} \mathcal{L} \sqrt{\tau}-1$ reads
$$\delta'=\dfrac{\sqrt{2}}{2}  \sqrt{\tau},\,\, \tilde{t}'=\dfrac{\sqrt{2}}{4}\mathcal{L}/ \sqrt{\tau},$$
with other variables conserved. Moreover the 0-energy hypersurface is given by the relationship $\dfrac{\sqrt{2}}{2}\mathcal{L} \sqrt{\tau}-1=0$, or equivalently $\tau=2 \mathcal{L}^{-2}$. The equation for $\tilde{t}$ thus becomes $\tilde{t}'=\dfrac{1}{4} \mathcal{L}^{2}$. Recall that the negative of the Keplerian energy satisfies $\tau=\dfrac{1}{2 a}$, where $a$ is the semi-major axis of the instantaneous Keplerian elliptic orbit. We therefore have $\mathcal{L}=2 \sqrt{a}$.

\subsection{The $\kappa$-Localization}

{In this section we  make the assumption that we always restrict the system to a fixed open set
$$\mathcal{U}:=\{0<\mathcal{L} < \tilde{L}_{\ast}, \tau > \tilde{\tau}_{\ast}\}$$ 
of the regularized extended phase space, for some fixed $L_{\ast}>0$ and $\tilde{\tau}_{\ast}>0$. The argument below shows that in \eqref{eq: eq motion perturbed Kepler}, when $a$ and $\tau^{-1}$ are small enough, any orbit of $H_{\varepsilon}$ in $H^{-1}_{\varepsilon}(0)$ starting with such semi major axis $a$ and $\tau$ will stay in this neighborhood for the fictitious time $s \in [0, S]$, in which $S$ is the fictitious time required for the $\tilde{t}$-variable to change from $0$ to $1$. All the norms of $U$ (or rather $\tilde{U}$) and its partial derivatives, which depend also on the external parameter $\varepsilon \in [0, 1]$ are taken as $L^{\infty}$ norms in $\mathcal{U} \times [0, 1]$. By properly shrinking the symplectic chart when necessary, we may assume that these norms all take bounded values.}

We set $\mathcal{L}(0)=\kappa$ for $\kappa>0$ sufficiently small. We shall explain in this subsection that any sufficiently small $\kappa$ gives a localized system which gives several a priori bounds on the corresponding periodic solutions of period $1$ in the $t$ variable. 



First, we have the energy constraint 
\begin{equation}\label{eq: energ const}
\dfrac{\sqrt{2}}{2} \kappa \sqrt{\tau(0)} + \kappa^{2} \tau^{-1}(0) \varepsilon \tilde{U}=1.
\end{equation}
When
$\lim_{\kappa \to 0} \sqrt{\tau(0)} \kappa = 0$, it follows from this energy constraint that $\tau(0) \sim \kappa^{2} \varepsilon U$ which is not allowed for small $\kappa$ by our definition of the neighborhood $\mathcal{U}$. 

We therefore necessarily have $\lim_{\kappa \to 0} \sqrt{\tau(0)} \kappa \neq 0$. We then deduce from \eqref{eq: energ const} that 
$$\lim_{\kappa \to 0} \sqrt{\tau(0)} \kappa=\sqrt{2},$$ therefore $\dfrac{1}{\sqrt{\tau(0)}}=\dfrac{\sqrt{2}}{2} \kappa+o(\kappa)$. 

We thus have $\kappa^{2} \tau^{-1}(0)=O(\kappa^{4})$. By plugging these into \eqref{eq: energ const}, we find that  actually $\dfrac{1}{\sqrt{\tau(0)}}=\dfrac{\sqrt{2}}{2} \kappa+O(\kappa^{3})$. Thus
$$\dfrac{\sqrt{2}}{\sqrt{ \tau(0)}} \in (\kappa-\kappa^{2}, \kappa+\kappa^{2})$$ 
for $\kappa$ sufficiently small. 

The region 
$$\Bigl\{\mathcal{L}(s), \dfrac{\sqrt{2}}{\sqrt{\tau(s)}} \in (\kappa-\kappa^{3/2}, \kappa+\kappa^{3/2})\Bigr\} $$
on $\{H_{\varepsilon}=0\}$ is referred to as the $\kappa$-localization. We shall show that

\begin{lemma}\label{prop: loc ansatz} When $\kappa$ is sufficiently small, any 0-energy solution of $H_{\varepsilon}$ with $\mathcal{L}(0)=\kappa$ lies entirely in the  $\kappa$-localization for $s \in [0, S]$ such that $\int_{0}^{S} \tilde{t}' d s=1$. 
\end{lemma}

We take these as a localization ansatz for the following estimates. 

To prove Lem. \ref{prop: loc ansatz}, we shall need the weaker inclusion 
$$\mathcal{L}(s), \dfrac{\sqrt{2}}{\sqrt{\tau(s)}} \in \Bigl(\dfrac{1}{2} \kappa, \dfrac{3}{2} \kappa\Bigr), s \in [0, S],$$
which we refer to as the weak localization ansatz. When the weak localization ansatz holds, we have that for sufficiently small $\kappa$ for all $\varepsilon \in [0, 1]$, that 
$$\mathcal{L}^{'}=-\mathcal{L}^{2} \tau^{-1} \varepsilon \Bigl(\dfrac{\partial \tilde{U}}{\partial \delta}+\dfrac{\partial \tilde{U}}{\partial t} \dfrac{\partial t}{\partial \delta}\Bigr).$$
Note that $\dfrac{\partial t}{\partial \delta}=\mathcal{L} \tau^{-1} \dfrac{\partial \tilde{\phi}}{\partial \delta}=O(\kappa^{3})$, therefore there exists $C_{1}>0$ such that
$$ |\mathcal{L}'| \le \dfrac{3^{4}}{2^{5}}\kappa^{4} \Bigl(\|\dfrac{\partial \tilde{U}}{\partial \delta} \|_{\infty}+ \dfrac{3^{3}}{2^{4}}\kappa^{3} \| \dfrac{\partial \tilde{U}}{\partial t}\|_{\infty} \| \dfrac{\partial \tilde{\phi}}{\partial \delta}\|_{\infty}\Bigr) \le C_{1} \kappa^{4}.$$
We also have
$$\delta'=\dfrac{\sqrt{2}}{2} \sqrt{\tau} + \mathcal{L} \tau^{-1} \varepsilon \bigl(2 \tilde{U} +\mathcal{L} \dfrac{\partial \tilde{U}}{\partial \mathcal{L}}\bigr).$$
Note that 
$$\dfrac{\partial \tilde{U}}{\partial \mathcal{L}}=\dfrac{\partial \hat{U}}{\partial t} \dfrac{\partial t}{\partial \mathcal{L}}+O(\tau^{-1/2})=O(\kappa), $$
thus for sufficiently small $\kappa$ there holds
$$|\delta'|  \le 3 \kappa^{-1}. $$

Next we consider $\hat{\xi}^{'}$ and $\hat{\zeta}'$. We have
$$\hat{\xi}^{'}=-\mathcal{L}^{2} \tau^{-1} \varepsilon \bigl(\dfrac{\partial \tilde{U}}{\partial \zeta}+\dfrac{\partial \tilde{U}}{\partial t} \dfrac{\partial t}{\partial \zeta}\bigr) \dfrac{\partial \zeta}{\partial \hat{\zeta}},$$
therefore
$$ |\hat{\xi}'| \le \dfrac{3^{4}}{2^{3}}\kappa^{7/2} \|\dfrac{\partial \tilde{U}}{\partial \zeta}+\dfrac{\partial \tilde{U}}{\partial t} \dfrac{\partial t}{\partial \zeta} \|_{\infty} := C_{2} \kappa^{7/2}, C_{2}>0;$$
similarly
$$\hat{\zeta}^{'}=\mathcal{L}^{2} \tau^{-1} \varepsilon \bigl(\dfrac{\partial \tilde{U}}{\partial \xi}+\dfrac{\partial \tilde{U}}{\partial t} \dfrac{\partial t}{\partial \xi}\bigr)\dfrac{\partial \xi}{\partial \hat{\xi}},$$
thus
$$|\hat{\xi}'| \le \dfrac{3^{4}}{2^{3}}\kappa^{7/2} \|\dfrac{\partial \tilde{U}}{\partial \xi}+\dfrac{\partial \tilde{U}}{\partial t} \dfrac{\partial t}{\partial \xi} \|_{\infty} := C_{3} \kappa^{7/2}, C_{3}>0.$$
Next we have
$$\tau'=-\mathcal{L}^{2} \tau^{-1} \varepsilon \dfrac{\partial \tilde{U}}{\partial t} \dfrac{\partial t}{\partial \tilde{t}} \hbox{  thus } |\tau'| \le \dfrac{3^{4}}{2^{4}}\kappa^{4} \bigl\|\dfrac{\partial \tilde{U}}{\partial \tilde{t}} \bigr\|_{\infty} := C_{4} \kappa^{4}, C_{4}>0.  $$
The right hand sides of the equations for $\mathcal{L'}$ and $\tau'$ are both of order $O(\kappa^{4})$. Therefore in view of our estimates for $\mathcal{L}(0) $ and $\tau(0)$,  we conclude that for sufficiently small $\kappa$, the weak localization ansatz holds for all $s \in [0, \kappa^{-3}]$.

Moreover, the equation for $\tilde{t}'$ takes the form:

\begin{eqnarray}
\tilde{t}'=\dfrac{\sqrt{2}}{4} \dfrac{\mathcal{L}}{\sqrt{\tau}}-\mathcal{L}^{2} \tau^{-1} \varepsilon (\tau^{-1}\tilde{U}-\dfrac{\partial \tilde{U}}{\partial \tau}).
\end{eqnarray}
Note that $\dfrac{\partial \tilde{U}}{\partial \tau}=\dfrac{\partial \tilde{U}}{\partial t} \dfrac{\partial t}{\partial \tau}$ is of order $\mathcal{L} \tau^{-2}$ by \eqref{eq: variable change LC}.  Therefore the last term in the right hand side of the above formula is of order $O(\kappa^{6})$. 
 
We therefore have 
$$|\tilde{t}'-\dfrac{\kappa^{2}}{4} | \le C_{5 }\kappa^{6}.$$


The time $S$ is given implicitly by the integral $\int_{0}^{S} \tilde{t}' ds=1$. 
It follows from the above estimate for $\tilde{t}'$ that
$$|S-4 \kappa^{-2}| \le  C_{6} \kappa^{2}, C_{6} >0.$$

In particular the weak localization ansatz is valid for all $s \in [0, S]$ for $\kappa$ sufficiently small. We may therefore estimate
$$\max_{s \in [0, S]}|\mathcal{L}(s)-\mathcal{L}(0)| \le \int_{0}^{S} |\mathcal{L}'| d s \le C_{1} \kappa^{4} \cdot 5 \kappa^{-2}=5 C_{1} \kappa^{2}, $$
and therefore for $\mathcal{L}(0)=\kappa$ it holds that 
$$\mathcal{L}(s) \in (\kappa-\kappa^{3/2},\kappa+\kappa^{3/2}), \qquad \hbox{ for } s \in [0, S]$$ for $\kappa$ sufficiently small. Likewise, for $\dfrac{\sqrt{2}}{\sqrt{\mathcal{\tau}(0)}} \in (\kappa-\kappa^{3}, \kappa+\kappa^{3})$ it holds that 
$$\dfrac{\sqrt{2}}{\sqrt{\mathcal{\tau}(s)}} \in (\kappa-\kappa^{3/2},  \kappa+\kappa^{3/2}), \qquad s \in [0, S]$$ 
for $\kappa$ sufficiently small. This proves Lem. \ref{prop: loc ansatz}. 

Moreover, we have 
$$|\hat{\xi}(0)|, |\hat{\zeta}(0)| = O(\sqrt{\mathcal{L}})=O(\kappa^{1/2})$$
and
$$\max_{s \in [0, S]}|\hat{\xi}(s)-\hat{\xi}(0)| \le \int_{0}^{S} |\hat{\xi}'| d s \le C_{2} \kappa^{7/2} \cdot 5 \kappa^{-2}=5 C_{2} \kappa^{3/2}, $$
$$\max_{s \in [0, S]}|\hat{\zeta}(s)-\hat{\zeta}(0)| \le \int_{0}^{S} |\hat{\zeta}'| d s \le C_{3} \kappa^{7/2} \cdot 5 \kappa^{-2}=5 C_{3} \kappa^{3/2}. $$
Thus for sufficiently small $\kappa>0$ we may indeed argue within a local symplectic chart as we do now.




The action functional along a solution curve (not necessarily periodic) from  $\tilde{t}=0$ to $\tilde{t}=1$ with energy zero in these variables is computed as 
$$\mathcal{A}:=\int_{0}^{S}  (\mathcal{L} \delta'+ \hat{\xi} \hat{\zeta}'+\tau \tilde{t}' )d s.$$ 

For the unperturbed system $H_{0}$, if we write $\kappa$ in place of $\mathcal{L}=\mathcal{L}_{n}=2^{2/3} \pi^{-1/3}  n^{-1/3}$, we obtain the order estimate $\mathcal{A}_{0} \sim \kappa^{-2}$.


With the estimates established above, we also have that 
\begin{lemma}\label{lem: 3.1}
When $\kappa>0$ is sufficiently small, it holds that
$\mathcal{A}_{\varepsilon} \sim \kappa^{-2}$
for any $\varepsilon \in [0, 1].$
\end{lemma}

{Note that due to the additional shift in the transformation $t \mapsto \tilde{t}$, in general it does not hold that 
$$t(S)-t(0)=1 \Leftrightarrow \tilde{t}(S)-\tilde{t}(0)=1.$$
Nevertheless it is enough to observe from \eqref{eq: variable change LC} that this equivalence holds for $S$-periodic solutions (in the time variable $s$) of the regularized system.}


We therefore draw the following conclusions for periodic orbits of the forced system with $t$-period $1$, which justify that our analysis is valid in any prescribed small neighborhood of the origin as long as $\kappa>0$ is sufficiently small. \begin{prop}\label{prop: control 0} For all $\varepsilon \in [0, 1]$, when $\kappa>0$ is sufficiently small, the action $\mathcal{A}_{\varepsilon}$ and $\tau$ as computed along a periodic solution from an initial value with $\mathcal{L}(0)=\kappa$ for the time interval $[0, S]$ such that 
$$t(0)=0, t(S)=1$$
are of the same order $\kappa^{-2}$. The quantity $|q(s)|$ is of the order $O(\kappa^{2})$ for $s \in [0, S]$, and are thus confined to a sufficiently small neighborhood of the origin. 
\end{prop}

In particular it is enough for our analysis to have that the function $U(q, t)$ or $U(q, t, \varepsilon)$ be regularly defined in a sufficiently small neighborhood of the origin.




\subsection{Extension to Moser-regularized Systems in All Dimensions}
To show the validity of Lem. \ref{prop: loc ansatz}, Lem. \ref{lem: 3.1} and Prop. \ref{prop: control 0} within the context of Moser regularization in any dimension, we observe that in dimension 2 this is obtained by simply going through the 2-to-1 cover. 

In higher dimensions, with the set of symplectic variables as constructed in Section \ref{Sec: 2.2.3}, we see that the estimates for the variables $(\mathcal{I}, \theta, \tau, \tilde{t})$ are valid in the same way in all dimensions, and the additional variables $(\hat{\xi}_{i}, \hat{\zeta}_{i})$ satisfies the same estimates as the pair $(\hat{\xi}, \hat{\zeta})$ in the Levi-Civita regularized system. With this argument we conclude that

\begin{prop}\label{prop: all d} Lem. \ref{prop: loc ansatz}, Lem. \ref{lem: 3.1} and Prop. \ref{prop: control 0} hold true for the Moser regularized system of  \eqref{eq: eq motion perturbed Kepler} in any dimensions.
\end{prop}

\section{Rescalings of Periodic Manifolds and Proof of Theorem \ref{thm: main}}\label{Sec: 4}








\subsection{The Rabinowitz Action Functional and Rescaling} 
Let $(X, \omega= d \lambda)$ be an exact symplectic manifold. Set
$$(\hat{X}, \hat{\omega}=d \hat{\lambda})=(X \times T^{*} S^{1}, \omega \oplus d \tau \wedge d t),$$
where we have denoted by $(\tau, t)$ the variables in $T^{*} S^{1}$. On $(\tilde{X}, \tilde{\omega})$ is defined a Hamiltonian $H \in C^{\infty} (\hat{X}, \R)$. By a periodic orbit we mean a closed orbit of $(\hat{X}, \hat{\omega}, H)$ with energy 0 along which the $t$-variable winds exactly once in the corresponding $S^{1}$-factor. 

Let $u \in C^{\infty} (S^{1}, \hat{X})$ be a smooth loop in $\hat{X}$, where $S^{1}=\R/ \Z$. In this subsection, we normalize the (fictitious) time along a closed orbit in the extended phase space to have period $1$ in this section, and denote it by $\underline{t}$. 
The (original, non-normalized) fictitious time in the extended phase space will not appear in this section.

For $(u, \eta) \in C^{\infty}  (S^{1}, \hat{X}) \times (0, \infty)$, the Rabinowitz action functional of $u$ with respect to a Hamiltonian function $H$ is defined as 
$$\mathcal{A}^{H} (u, \eta)=-\int_{0}^{1} u^{*} \hat{\lambda} +\eta \int_{0}^{1} H(u) d t, $$
with critical point equations
$$\dfrac{d u}{d \underline{t}}=\eta X_{H}(u), \qquad \int_{0}^{1} H(u) d t =0 \Leftrightarrow H(u) \equiv 0.$$
In other words, the critical points of the Rabinowitz action functional are 0-energy periodic orbits of $H$ with period $\eta$. 


A critical manifold, \emph{i.e.} a manifold consisting of critical points of the action functional, is called \emph{Morse-Bott} if at each point the kernel of the Hessian of the action functional  agrees with the tangent space at the point. An explicit computation of the kernel of the Hessian of the Rabinowitz action functional in \cite[Chapter 7, Section 3]{Frauenfelder-Koert} directly shows that 

\begin{prop} A non-degenerate periodic manifold in the sense of Condition \eqref{Condition: Nondegeneracy} is exactly a Morse-Bott critical manifold of the associated Rabinowitz action functional. 
\end{prop}

Therefore, in any dimensions, the periodic manifolds $\Lambda_{n}$ of $H_{0}$ are Morse-Bott critical manifolds.



We now investigate how the Rabinowitz action functional and its critical point equations change under rescaling. We take a conformal symplectic mapping
$$\phi_{\kappa}: \hat{X} \to \hat{X}$$
such that
$$\tilde{\lambda}=\kappa^{2} \phi_{\kappa}^{*} \hat{\lambda}, \qquad \tilde{\omega}=\kappa^{2} \phi_{\kappa}^{*} \hat{\omega}$$
for $\kappa>0$.

We set 
$$\tilde{u}:=\phi_{\kappa}^{-1} \circ u \in C^{\infty}(S^{1}, \hat{X}).$$
We 
set $\tilde{H}_{\kappa}(\tilde{u}):
=\phi_{\kappa}^{*} H(u)$.
Accordingly we define $\eta=\kappa^{-2} \tilde{\eta}$.

{We compute} 
$$\int_{0}^{1} u^{*} \hat{\lambda}=\int_{0}^{1} (\phi_{\kappa} \circ \tilde{u})^{*} (\kappa^{-2} {\phi_{\kappa}^{-1}}^{*}\tilde{\lambda})=\kappa^{-2} \int_{0}^{1} \tilde{u}^{*} \tilde{\lambda}.$$

Thus we have that
$$\mathcal{A}^{H} (u, \eta)=-\kappa^{-2} \int_{0}^{1} \tilde{u}^{*} \tilde{\lambda} + \kappa^{-2} \tilde{\eta}  \int_{0}^{1} 
\tilde{H}_{\kappa}(\tilde{u}) d t,$$
or equivalently
$\mathcal{A}^{H}(u, \eta)=\kappa^{2} \mathcal{A}^{\tilde{H}_{\kappa}} (\tilde{u}, \tilde{\eta})$
with critical point equations
$$\dfrac{d}{d \underline{t}} \tilde{u}=\tilde{\eta} X_{\tilde{H}_{\kappa}}(\tilde{u}), \qquad \int_{0}^{1} \tilde{H}_{\kappa}(\tilde{u}) d t =0 \Leftrightarrow \tilde{H}_{\kappa}(\tilde{u}) \equiv 0.$$

Therefore $\tilde{u}$ is a critical point of the rescaled Rabinowitz action functional, with the rescaled multiplier $\tilde{\eta}$.




\subsection{Rescaled $\kappa$-Localization}
In our problem we have
$$\hat{X}=T^{*} S^{d} \times T^{*} S^{1}.$$ 
We take a $\kappa$ localization  for $\kappa=\kappa_{n}=2^{2/3} \pi^{-1/3} n^{-1/3}$ sufficiently small in the phase space. Precisely this means that on the zero-energy hypersurface there holds
$$\mathcal{L}, \dfrac{\sqrt{2}}{\sqrt{\tau}} \in (\kappa_{n}-\kappa_{n}^{3/2}, \kappa_{n}+\kappa_{n}^{3/2}).$$
We make a conformal change of the symplectic variables as
\begin{equation}\label{eq: change of coordinates}
 (\mathcal{L} =\kappa_{n} \tilde{\mathcal{L}}, \mathcal{\delta} = \kappa_{n}^{-3} \tilde{\delta}, \hat{\xi}_{i}=\kappa_{n}^{-1} \tilde{\xi}_{i},  \hat{\zeta}_{i}=\kappa_{n}^{-1} \tilde{\zeta}_{i},  \tau=\kappa_{n}^{-2} \tilde{\tau}, \tilde{t}=\tilde{t}).
 \end{equation}

We have 
$$\tilde{\lambda}=( \tilde{\mathcal{L}}  d \tilde{\delta}+\sum_{i=1}^{d-1} \tilde{\xi}_{i} d \tilde{\zeta}_{i} +  \tilde{\tau}  d \tilde{t})=\kappa_{n}^{2}(\mathcal{L}  d \delta+\sum_{i=1}^{d-1} \tilde{\xi}_{i} d \tilde{\zeta}_{i}+  \tau d \tilde{t})=\kappa_{n}^{2} \hat{\lambda},$$
thus
$$\tilde{\omega}=d \tilde{\lambda}=\kappa_{n}^{2} d \hat{\lambda}=\kappa_{n}^{2} \hat{\omega.}$$

In other words, the mapping 
$$\phi_{\kappa_{n}}: \hat{X} \to \hat{X},  (\tilde{L}, \tilde{\delta}, \tilde{\xi}, \tilde{\zeta}, \tilde{\tau}, \tilde{t}) \mapsto (L, \delta,  \hat{\xi}, \hat{\zeta},  \tau, \tilde{t})$$ 
satisfies
$$\tilde{\lambda}=\kappa_{n}^{2} \phi_{\kappa_{n}}^{*} \hat{\lambda}, \qquad \tilde{\omega}=\kappa_{n}^{2} \phi_{\kappa_{n}}^{*} \hat{\omega}.$$


Recall that we have $H_{0}=\dfrac{\sqrt{2}}{2} \mathcal{L} \sqrt{\tau}-1$. 
It is instructive to write 
$$H_{\varepsilon}=H_{0}+O(\kappa_{n}^{4};\varepsilon),$$
as it follows from \eqref{eq: estimate on U} that the perturbation term is of order $O(\kappa_{n}^{4})$ in $\kappa_{n}$ and tends to zero when $\varepsilon \to 0$. After rescaling, we have
that
\begin{equation}\label{eq: normalized kappa}
\tilde{H}_{\varepsilon, \kappa_{n}}=\dfrac{\sqrt{2}}{2} \tilde{\mathcal{L}} \sqrt{\tilde{\tau}}-1+O(\kappa_{n}^{4};\varepsilon).
\end{equation}

 
We now apply this argument for each $n$ sufficiently large to normalize the corresponding Keplerian periodic manifolds $\Lambda_{n}$ of $H_{0}$, which we describe in local charts, to 
$$\Lambda_{1}^{(n)}:=\{\tilde{\mathcal{L}}=\mathcal{L}_{1}, \tau=\tau_{1}, \tilde{\delta} \in \R/ (2 \pi/n) \Z\}.$$
{Globally $\Lambda_{1}^{(n)}$ is realized in the fibre bundle $B^{(n)}$ with fibre $T_{+}^{*} (\R/ (2 \pi/n) \Z \times S^{1})$ over the normalized orbit space $\Omega_{\mathcal{O}}$ as
$$\Lambda_{1}^{(n)}:=\{\tilde{\mathcal{L}}=\mathcal{L}_{1}, \tau=\tau_{1}\} \subset B^{(n)}.$$}

{To further uniformize the periodic manifolds we now lift the systems $H_{0}$ and $H_{\varepsilon}$ defined on $B^{(n)}$ to a fibrewise n-cover of $B^{(n)}$ by lifting $T_{+}^{*} (\R/ (2 \pi/n) \Z \times S^{1})$ to $T_{+}^{*} \T^{2}$  fibrewisely. The fibres $T_{+}^{*} (\R/ (2 \pi/n) \Z \times S^{1})$ and $T_{+}^{*} \T^{2}$ are both quotients of $T_{+}^{*} \R^{2}$ by $\Z^{2}$-actions on $\R^{2}$, which can be identified by means of an affine transformation on $\R^{2}$. Therefore both $B^{(n)}$ and its fibrewise n-cover are quotients of the same $T_{+}^{*} \T^{2}$-bundle with respect to the (different but identifiable) actions of $\Z^{2}$ which act fibrewisely. They are therefore isomorphic as fibre bundles. On the other hand, all $B^{(n)}$'s are isomorphic. Therefore for all $n$, the fibrewise n-covers of $B^{(n)}$ are isomorphic to the fibre bundle $B:=B^{(1)}$ with fibres  $T^{*}_{+} \T^{2}$ over the base  $\Omega_{\mathcal{O}}$. Moreover, the equivariant lift of the symplectic structure on $B^{(n)}$ agrees with the canonical symplectic structure on $B^{(1)}$. We may therefore take $B:=B^{(1)}$ with its canonical symplectic structure as a uniform representative. The periodic manifolds $\Lambda_{1}^{(n)}$ are now lifted to $\Lambda_{1}=\Lambda_{1}^{(1)}$ in $B$.}


 
The unperturbed function $H_{0}$ is independent of the angles and can be lifted directly.
{The perturbation term can be written in a chart as a function of the rescaled variables
$$(\tilde{\mathcal{L}}, n \tilde{\delta}, \hat{\xi}_{i}, \hat{\zeta}_{i}, \tilde{\tau}, \tilde{t} ),$$
in which in particular the angle $\tilde{\delta}$ is defined in $\R/ (2 \pi/n) \Z$.
which then lifts to a function depending on the variables
$$(\tilde{\mathcal{L}}, [n \tilde{\delta}], \hat{\xi}_{i}, \hat{\zeta}_{i}, \tilde{\tau}, \tilde{t} ),$$
while in the latter 
$$[n \tilde{\delta}]=n \tilde{\delta} \mod 2 \pi$$ 
and the angle $\tilde{\delta}$ is extended to be defined in $\R/ (2 \pi) \Z$.}

{In an overlap region of two charts given respectively by
$$(\tilde{\mathcal{L}}, \tilde{\delta}_{1}, \tilde{\xi}_{i, 1}, \tilde{\zeta}_{i, 1}, \tilde{\tau}, \tilde{t}_{1} ),$$
and
$$(\tilde{\mathcal{L}}, \tilde{\delta}_{2}, \tilde{\xi}_{i, 2}, \tilde{\zeta}_{i, 2}, \tilde{\tau}, \tilde{t}_{2} ),$$
the transition map is induced from 
the corresponding transition map of the charts given by the non-rescaled variables
$$\bar{\phi}: ({\mathcal{L}}, {\delta}_{1}, \hat{\xi}_{i, 1}, \hat{\zeta}_{i, 1}, {\tau}, \tilde{t}_{1} ) \mapsto ({\mathcal{L}}, {\delta}_{2}, \hat{\xi}_{i, 2}, \hat{\zeta}_{i, 2}, {\tau}, \tilde{t}_{2})$$ 
which in particular implies that
$$({\mathcal{L}}, n \tilde{\delta}_{2}, \hat{\xi}_{i, 2}, \hat{\zeta}_{i, 2}, {\tau}, \tilde{t}_{2})=\bar{\phi}  ({\mathcal{L}}, n \tilde{\delta}_{1}, \hat{\xi}_{i, 1}, \hat{\zeta}_{i, 1}, {\tau}, \tilde{t}_{1} )$$
in which $\tilde{\delta}_{1}, \tilde{\delta}_{2} \in \R/ (2 \pi/n) \Z$. 
This equality extends to
$$({\mathcal{L}}, [n \tilde{\delta}_{2}], \hat{\xi}_{i, 2}, \hat{\zeta}_{i, 2}, {\tau}, \tilde{t}_{2})=\bar{\phi}  ({\mathcal{L}}, [n \tilde{\delta}_{1}], \hat{\xi}_{i, 1}, \hat{\zeta}_{i, 1}, {\tau}, \tilde{t}_{1} )$$
in which we now regard $\tilde{\delta}_{1}, \tilde{\delta}_{2} \in \R/ 2 \pi \Z$. 

This shows that the lift of the perturbation term takes consistent values on the overlap region of two lifted charts and is thus globally defined in a neighborhood of $\Lambda_{1}$ in $B$. Moreover it still assumes the form $O(\kappa_{n}^{4};\varepsilon)$. In particular, this lift construction does not change its $C^{1}$-norm.}


We now remark on the $C^{1}$-norm of the perturbation term, which assumes the form $O(\kappa_{n}^{4};\varepsilon)$. We see that due to the rescaling of the variables \eqref{eq: change of coordinates}, in particular the fast angle $\delta$, the partial derivative of the perturbation with respect to $\tilde{\delta}$ is only of order $O(\kappa_{n})$ while other partial derivatives remain to be of order $O(\kappa_{n}^{4})$. As a result, the $C^{1}$-norm of the pertubation term is of order $O(\kappa_{n})$. In general, we may only have an estimate of the $C^{2}$-norm of the perturbation as $O(\kappa_{n}^{-2})$. This estimate is not bounded when $\kappa_{n} \to 0$. Since the theorem of Weinstein \cite{Weinstein} requires a proper smallness of the $C^{2}$-norm of the perturbation, this can only be applied to obtain our result in special cases. To obtain the result in full generality we shall instead conclude with the local Rabinowitz-Floer homology argument from Appendix \ref{Appendix A}. 

\smallskip

\subsection{Proof of  Theorem \ref{thm: main} and Corollary \ref{Cor: main}}

\textbf{Proof of Theorem \ref{thm: main}}
We apply Thm. \ref{thm: A} to the rescaled system \eqref{eq: normalized kappa} near the normalized periodic manifold $\Lambda_{1}$ of  \eqref{eq: normalized kappa} with $\varepsilon=0$. The period of periodic orbits from $\Lambda_{1}$ is $S_{1}=(2 \pi)^{2/3}$. We take $T^{\pm}=S_{1} \pm 10^{-4}$. Moreover, we have $K^{+}, K^{-} =O(\kappa_{n}^{4})$. In this setting Thm. \ref{thm: A} is applicable for all sufficiently small $\kappa_{n}$, which then assures the existence of a periodic orbit with action $\mathcal{A}_{0}^{(1)}+O(\kappa_{n}^{4})$ of the system \eqref{eq: normalized kappa} for all sufficiently large $n$. After rescaling back this corresponds to a periodic orbit of the system \eqref{eq: eq motion perturbed Kepler} with action of the order $\kappa_{n}^{-2}+O(\kappa_{n}^{2})$. Viewing from Prop. \ref{prop: action value} that the action gaps of the regularized Kepler problem in extended phase space $H_{0}$ is of order $\tilde \kappa_{n}$ for sufficiently small $\kappa_{n}$, the periodic orbits thus obtained take infinitely many distinct action values. Since the perturbations to the Moser-regularized Kepler flow vanish at the collisions,  any of these orbits passes transversally through the set of collisions along which the variable $\tau$ is continuous, we conclude that this orbit gives rise to a generalized periodic orbit of the system \eqref{eq: eq motion perturbed Kepler} resp. \eqref{eq: eq motion perturbed Kepler varepsilon=1} with the same action value. Consequently we obtain infinitely many generalized periodic orbits of the system \eqref{eq: eq motion perturbed Kepler} resp. \eqref{eq: eq motion perturbed Kepler varepsilon=1}.
\hfill $\square$

Indeed the argument shows that there exists $N^{*}>0$ which depends only on the $C^{1}$-norm of $U$, such that for all $n >N^{*}$ there bifurcates at least one periodic orbit  from each periodic manifold $\Lambda_{n}$. For $n \le N^{*}$ we may control the smallness of the $C^{1}$-norm of the perturbation by choosing small $\varepsilon>0$ as in \cite{Boscaggin--Ortega--Zhao}. We thus obtain the following corollary:

\begin{cor}\label{cor: main 1}There exists $\varepsilon^{*}>0$ depending only on the $C^{1}$-norm of $U$, such that for $\varepsilon \in [0, \varepsilon^{*}]$ the system \eqref{eq: eq motion perturbed Kepler} has infinitely many periodic orbits bifurcating from each of the periodic manifolds $\{\Lambda_{n}\}_{n=1,2,\cdots}$ . 
\end{cor}

\begin{rem} We have seen that in general the $C^{2}$-norm of $\varepsilon U$ is of the order $O(\kappa^{-2};\varepsilon) $, caused by the rescaling of the fast angle $\delta$. This can be improved provided that the perturbation is independent of $\delta$, or when the perturbation is of the order $o(\kappa^{6}; \varepsilon)$. The latter is achieved, for example, when the function $U(q, t, \varepsilon)$ is super-cubic in $|q|$, \emph{i.e.} $U(q, t, \varepsilon)=O(|q|^{4}; \varepsilon)$ in a neighborhood of the origin. In these cases we may as well conclude Thm. \ref{thm: main} and Cor. \ref{cor: main 1} with Weinstein's theorem. 
\end{rem}

{Finally we remark that by realizing circular or elliptic restricted three-body problems as periodically forced Kepler problems  as in Appendix \ref{Appendix C}, we obtain infinitely many generalized periodic orbits accumulating each of the primaries, which proves Cor. \ref{Cor: main}}.

\medskip
\medskip
{\bf Acknowledgements} Many thanks to Urs Frauenfelder for essential contribution to Appendix A and for many helpful and stimulating discussions, to Rafael Ortega for helpful discussions and to Alain Chenciner for various corrections to previous versions of this manuscript. The author is supported by DFG ZH 605/1-1.

\appendix

\newpage

\section{$C^{1}$-Persistence of Periodic Orbits via a Localized Homotopy-Streching Argument}\label{Appendix A}
\qquad \qquad \qquad \qquad \qquad \qquad {Urs Frauenfelder\footnote{University of Augsburg, Augsburg, DE: urs.frauenfelder@math.uni-augsburg.de}, Lei Zhao} 
\medskip

The theory of (periodic) Floer Homology provides powerful tools to detect periodic orbits in a Hamiltonian system. The Floer theory associated to the Rabinowitz action functional is particularly helpful to detect periodic orbits on a given energy hypersurface of an autonomous Hamiltonian system. Precisely speaking for the Rabinowitz Floer Homology to be well-defined one typically requires that the energy hypersurface is compact, which does not hold in our case. In this Appendix we use a localized Rabinowitz-Floer Homology argument to show the existence of a periodic orbit bifurcating from a Morse-Bott periodic manifold under $C^{1}$-small smooth perturbations of the Hamiltonian. 

We assume familiarity with the theory of Floer Homology which can be found in standard references such as \cite{McDuff-Salamon}, \cite{Salamon}, \cite{AudinDamian}. Local versions of Floer homology has been previously used in \cite{GinzburgGurel}, \cite{HryniewiczMacarini}, \cite{GinzburgHeinHryniewiczMacarini}, \cite{McLean}, \cite{KimKimKwon}. Our argument in particular treats the Morse-Bott situation in Rabinowitz-Floer homology. 

Suppose that $(M,\omega)$ is a symplectic manifold and $H \colon M \to \mathbb{R}$ is a smooth
function referred to as the Hamiltonian. We do not require that $M$ is closed (\emph{i.e.} compact without boundary) but to simplify the discussion we suppose that the symplectic manifold is exact, i.e., $\omega=d\lambda$ for a 
one-form $\lambda$. We assume that $0$ is a regular value, such that $\Sigma=H^{-1}(0)$ is
a (not necessarily closed) codimension one submanifold of $M$. The choice of the regular value $0$ is certainly not restrictive, as we may always shift the Hamiltonian by an additive constant. The Hamiltonian vector field of $H$ is implicitly defined by the condition
$$dH=\omega(\cdot, X_H).$$
By our assumption $H$ is autonomous, i.e., independent of time, therefore $H$ remains constant
along integral curves of its Hamiltonian vector field $X_H$. 

Periodic orbits of $X_H$ on $H^{-1}(0)$ can be detected variationally as the critical points of the (Rabinowitz) action functional
$$\mathcal{A}^H \colon C^\infty(S^1,M) \times (0,\infty) \to \R, \quad
(v,\tau) \mapsto -\int v^*\lambda+\tau \int_0^1H(v(t))dt,$$
where $S^1=\mathbb{R}/\mathbb{Z}$. We suppose that
$$\mathcal{C} \subset C^\infty(S^1,M) \times (0,\infty)$$
is a closed (\emph{i.e.} compact without boundary) connected submanifold and is a Morse-Bott component of the critical set of $\mathcal{A}^H$, i.e.,
$$\mathcal{C} \subset \mathrm{crit}(\mathcal{A}^H), \qquad T_w \mathcal{C}=\ker \mathcal{H}_{\mathcal{A}^H}(w),\quad \forall\,\,w \in \mathcal{C},$$
where $\mathcal{H}_{\mathcal{A}^H}(w)$ denotes the Hessian of the action functional $\mathcal{A}^H$
at its critical point $w$. In particular
 $\mathcal{C}$ is isolated in the critical set $\mathrm{crit}(\mathcal{A}^H)$. 
Moreover, since the action is constant along $\mathcal{C}$, we may write
$$\mathcal{A}^H(\mathcal{C}):=\mathcal{A}^H(w), \quad w \in \mathcal{C}.$$
Since
$\mathcal{C}$ is closed, the set
$$M_\mathcal{C}=\big\{v(t):w=(v,\tau) \in \mathcal{C},\,\, t\in S^1\big\} \subset M$$
is compact and there exist finite
$$\tau_-:=\min\big\{\tau: w=(v,\tau) \in \mathcal{C}\big\},\qquad
\tau_+:=\max\big\{\tau: w=(v,\tau) \in \mathcal{C}\big\}.$$
We choose an open neighbourhood $V$ of $M_\mathcal{C}$ in $M$ with the property that its
closure $\overline{V} \subset M$ is compact. Furthermore, we choose positive real numbers
$T_-$ and $T_+$ such that
$$T_-<\tau_- \leq \tau_+<T_+.$$
We further use the following notation: If $K \in C^\infty(M,\mathbb{R})$ we introduce the continuous,
nonnegative functions
$$K^+ \colon M \to [0,\infty), \quad v \mapsto \max\big\{K(v),0\big\}$$
and 
$$K^- \colon M \to [0,\infty), \quad v \mapsto \max\big\{-K(v),0\big\}.$$ 
\begin{thm}\label{thm: A} \emph{There exists a $C^1$-open neighbourhood 
$\mathcal{U}=\mathcal{U}(V,T_-,T_+)$ of $0$ in $C^\infty(M,\mathbb{R})$ with the following property: 
For every $K \in \mathcal{U}$ 
 there exists a critical point $w=(v,\tau)$ of 
 $\mathcal{A}^{H+K}$ such that
\begin{eqnarray*}
& &\mathcal{A}^H(\mathcal{C})-T_+\big(\max K^+|_{\overline{V}}\,+\max K^-|_{\overline{V}}\big)\\
&\leq& \mathcal{A}^{H+K}(w)\\
&\leq& \mathcal{A}^H(\mathcal{C})+T_+\big(\max K^+|_{\overline{V}}\,+\max K^-|_{\overline{V}}\big)
\end{eqnarray*}
and
$$v(t) \in V,\quad t \in S^1, \qquad T_-<\tau<T_+.$$
 }
\end{thm}

We prove the Theorem by a homotopy-stretching argument.\footnote{Note that the variables $s, t$ in this Appendix do not refer to time variables in the main part of the article.} As a preparation for that purpose we first pick a bump function
$\beta \in C^\infty(\mathbb{R},[0,1])$ satisfying
$$\beta(s)=1,\quad s \in [-1,1],\qquad \beta(s)=0,\quad |s| \geq 2,$$
as well as
$$\beta'(s) \geq 0,\quad s \in [-2,-1],
\qquad \beta'(s) \leq0, \quad s \in [1,2].$$
We further choose a smooth family of compactly supported bump functions $\beta_r \in C^\infty(\mathbb{R},[0,1])$ for
$r \in [0,\infty)$ such that $\beta_0$ vanishes identically  and for $r \geq 1$ we have
$$\beta_r(s)=\left\{\begin{array}{cc}
\beta(s+r-1) & s \leq -r\\
1 & -r \leq s \leq r\\
\beta(s-r+1) & s \geq r.
\end{array}\right.$$ 
We further require that for every $r \in [0,\infty)$
$$\beta_r'(s) \geq 0, \quad s \leq 0, \qquad \beta_r'(s) \leq 0, \quad s\geq 0.$$
In particular, we have $\beta_1=\beta$. We fix $K \in C^\infty(M,\mathbb{R})$ and define
for $r \in [0,\infty)$ the time-dependent functional
$$\mathcal{A}_r=\mathcal{A}^{H,K}_r \colon C^\infty(S^1,M) \times (0,\infty) \times \mathbb{R}
\to \mathbb{R}$$
which maps $(v,\tau,s) \in C^\infty(S^1,M) \times (0,\infty) \times \mathbb{R}$ to 
\begin{eqnarray*}
\mathcal{A}_r(v,\tau,s)&:=&\mathcal{A}^H(v,\tau)+\beta_r(s) \tau \int_0^1 K(v(t))dt\\
&=&-\int v^*\lambda +\tau \int_0^1 \big(H+\beta_r(s)K\big)(v(t))dt.
\end{eqnarray*}
We choose an $\omega$-compatible almost complex structure $J$ on $M$, i.e., an almost complex structure
such that $\omega(\cdot,J\cdot)$ is a Riemannian metric on $M$. With the help of $J$ we endow the loop space
$C^\infty(S^1,M)$ with the $L^2$-metric obtained by integrating the Riemannian metric $\omega(\cdot,J\cdot)$
on $M$. On $C^\infty(S^1,M) \times (0,\infty)$ we consider the product of the $L^2$-metric on the loop space and
the standard metric on $(0,\infty)$. The (time-dependent) gradient of the time-dependent functional $\mathcal{A}_r$
with respect to this metric at a point $(v,\tau) \in C^\infty(S^1,M) \times (0,\infty)$  becomes
$$\nabla\mathcal{A}_r(v,\tau)(s)=\left(\begin{array}{cc}
J(v)(\partial_t v-\tau X_{H+\beta_r(s)K}(v))\\
\int_0^1 \big(H+\beta_r(s)K\big)(v)dt
\end{array}\right).
$$ 
We are interested in the maps
$$w=(v,\tau) \colon \mathbb{R} \to C^\infty(S^1, V) \times (T_-,T_+)$$
satisfying the gradient flow equations
\begin{equation}\label{grad}
\partial_s w(s)+\nabla\mathcal{A}_r(w(s))(s)=0, \quad s\in \mathbb{R}
\end{equation}
for $r \in [0,\infty)$ and the asymptotic conditions
\begin{equation}\label{asymp}
\lim_{s \to \pm \infty} w(s) \in \mathcal{C}.
\end{equation}
If we alternatively interpret these maps as
$$w=(v,\tau) \in C^\infty\big(\mathbb{R} \times S^1, V\big) \times C^\infty\big(\mathbb{R},(T_-,T+)\big),$$
then the gradient flow equation becomes 
\begin{equation}\label{grad2}
\left\{\begin{array}{c}
\partial_s v(s,t)+J\big(v(s,t)\big)\Big(\partial_t v(s,t)-\tau(s) X_{H+\beta_r(s)K}\big(v(s,t)\big)\Big)=0\\
\partial_s \tau(s)+\int_0^1 \big(H+\beta_r(s)K\big)\big(v(s,t)\big)dt=0.
\end{array}\right.
\end{equation}
\begin{lemma}\label{actest}
Suppose that $w=(v,\tau) \colon \mathbb{R} \to C^\infty(S^1, V) \times (T_-,T_+)$ satisfies the gradient flow 
equation (\ref{grad}) for $r \in [0,\infty)$ and the asymptotic conditions (\ref{asymp}). Then for every
$\sigma \in \mathbb{R}$ we have the estimate
\begin{eqnarray*}
& &\mathcal{A}^H(\mathcal{C})-T_+\big(\max K^+|_{\overline{V}}\,+\max K^-|_{\overline{V}}\big)\\
&\leq& \mathcal{A}_r\big(w(\sigma),\sigma\big)\\
&\leq& \mathcal{A}^H(\mathcal{C})+T_+\big(\max K^+|_{\overline{V}}\,+\max K^-|_{\overline{V}}\big).
\end{eqnarray*}
 
\end{lemma}
\textbf{Proof: } For $\sigma \in \mathbb{R}$, we estimate using the gradient flow equation and the asymptotic condition
\begin{eqnarray*}
\mathcal{A}_r\big(w(\sigma),\sigma\big)-\mathcal{A}^H(\mathcal{C})&=&\mathcal{A}_r(w(\sigma),\sigma)-\lim_{s \to -\infty} 
\mathcal{A}_r\big(w(s),s\big)\\
&=&\int_{-\infty}^\sigma \frac{d}{ds}\mathcal{A}\big(w(s),s\big)ds\\ 
&=&\int_{-\infty}^\sigma d \mathcal{A}_r\big(w(s),s\big)\partial_s w(s) ds+\int_{-\infty}^\sigma \partial_s \mathcal{A}_r\big(w(s),s\big)ds\\ 
&=&-\int_{-\infty}^\sigma \big|\big|\nabla \mathcal{A}_r\big(w(s),s)\big)\big|\big|^2 ds+\int_{-\infty}^\sigma \partial_s \mathcal{A}_r\big(w(s),s\big)ds\\ 
&\leq &\int_{-\infty}^\sigma \partial_s \mathcal{A}_r\big(w(s),s\big)ds.
\end{eqnarray*}
This implies
\begin{equation}\label{e1}
\mathcal{A}_r\big(w(\sigma),\sigma\big)\leq \mathcal{A}^H(\mathcal{C})+\int_{-\infty}^\sigma \partial_s \mathcal{A}_r\big(w(s),s\big)ds.
\end{equation}
To estimate the second term, we note that
$$\partial_s \mathcal{A}_r\big(w(s),s\big)=\beta'_r(s) \tau(s)\int_0^1 K(v(t,s))dt.$$
Note that $\tau(s)$ is always positive and is estimated from above by $T_+$. For $s \leq 0$ the derivative
of $\beta_r(s)$ is positive as well by our choice of the family of cutoff functions. 
Therefore we have
$$\partial_s \mathcal{A}_r\big(w(s),s\big) \leq \beta'_r(s) \cdot T_+ \cdot \max K^+|_{\overline{V}}, \qquad
s \leq 0.$$
For $s \geq 0$, we have by construction $\beta_r'(s) \leq 0$ and consequently
$$\partial_s \mathcal{A}_r\big(w(s),s\big) \leq -\beta'_r(s) \cdot T_+ \cdot \max K^-|_{\overline{V}}, \qquad
s \geq 0.$$
Combining these we obtain the estimate
\begin{eqnarray}\label{e2}
\int_{-\infty}^\sigma \partial_s \mathcal{A}_r\big(w(s),s\big)ds&\leq& \int_{-\infty}^{\min\{0,\sigma\}}
\beta'_r(s) \cdot T_+ \cdot \max K^+|_{\overline{V}}\,ds\\ \nonumber
& &-\int_0^{\max\{0,\sigma\}}\beta'_r(s)\cdot T_+ \cdot \max K^-|_{\overline{V}}\,ds\\ \nonumber
&\leq&T_+ \cdot \max K^+|_{\overline{V}}\int_{-\infty}^0
\beta'_r(s)ds\\ \nonumber
& &-T_+ \cdot \max K^-|_{\overline{V}}\int_0^\infty\beta'_r(s)ds\\ \nonumber
&=&T_+ \cdot \max K^+|_{\overline{V}} \cdot \beta_r(0)+T_+ \cdot \max K^-|_{\overline{V}}\cdot \beta_r(0)
\\ \nonumber
&\leq&T_+\big(\max K^+|_{\overline{V}}\,+\max K^-|_{\overline{V}}\big).
\end{eqnarray}
Plugging (\ref{e2}) into (\ref{e1}) we obtain the second inequality of the Lemma.
\\ 
To prove the first inequality we argue similarly while using the asymptotic at $\infty$ 
instead of $-\infty$: 
\begin{eqnarray*}
\mathcal{A}^H(\mathcal{C})-\mathcal{A}_r\big(w(\sigma),\sigma\big)&=&\lim_{s \to \infty} 
\mathcal{A}_r\big(w(s),s\big)-\mathcal{A}_r(w(\sigma),\sigma)\\
&=&-\int_\sigma^\infty \frac{d}{ds}\mathcal{A}\big(w(s),s\big)ds\\ 
&=&-\int_\sigma^\infty d \mathcal{A}_r\big(w(s),s\big)\partial_s w(s) ds-\int_\sigma^\infty \partial_s \mathcal{A}_r\big(w(s),s\big)ds\\ 
&=&\int_\sigma^\infty \big|\big|\nabla \mathcal{A}_r\big(w(s),s)\big)\big|\big|^2 ds-\int_\sigma^\infty \partial_s \mathcal{A}_r\big(w(s),s\big)ds\\ 
&\geq &-\int_\sigma^\infty \partial_s \mathcal{A}_r\big(w(s),s\big)ds.
\end{eqnarray*}
Therefore
\begin{equation}\label{e3}
\mathcal{A}_r\big(w(\sigma),\sigma\big)\geq \mathcal{A}^H(\mathcal{C})-\int_\sigma^\infty \partial_s \mathcal{A}_r\big(w(s),s\big)ds.
\end{equation}
The same reasoning as in the proof of (\ref{e2}) shows that
$$\int_\sigma^\infty \partial_s \mathcal{A}_r\big(w(s),s\big)ds \leq T_+\big(\max K^+|_{\overline{V}}\,+\max K^-|_{\overline{V}}\big),$$
so that the first inequality is proved as well. \hfill $\square$

We further choose 
$S_+, S_->0$ satisfying
$$T_-<S_-<\tau_- \leq \tau_+<S_+<T_+$$
and an open neighbourhood
$W$ of $M_\mathcal{C}$ with the property that
$$\overline{W} \subset V$$
and moreover, that every periodic orbit of $H$ of period less than or equal to $T_+$ contained in $\overline{W}$
belongs to $\mathcal{C}$. That such an open neighbourhood exists follows from the following reasoning: Otherwise
there exists a sequence of periodic orbits of period bounded by $T_+$ not belonging to $\mathcal{C}$ but converging
to $M_\mathcal{C}$. Because its period is uniformly bounded, by the Theorem of Arzela-Ascoli they have a convergent
subsequence which converges to a periodic orbit on $M_\mathcal{C}$. 
Therefore this limit orbit actually lies
in $\mathcal{C}$. But this contradicts the assumption that $\mathcal{C}$ is Morse-Bott. 
\begin{prop}\label{neighbour}
There exists a $C^1$-open neighbourhood $\mathcal{U}$ of $0$ in $C^\infty(M,\mathbb{R})$ with the following property:
Suppose that $w \colon \mathbb{R} \to C^\infty(S^1,V) \times (T_-,T_+)$ satisfies the gradient
flow equation (\ref{grad}) for $\mathcal{A}_r=\mathcal{A}_r^{H,K}$ with $r \in [0,\infty)$ and
$K \in \mathcal{U}$, as well as the asymptotic conditions (\ref{asymp}). Then
$$w(s) \in C^\infty(S^1,W) \times (S_-,S_+), \quad \forall\,\,s \in \mathbb{R}.$$
\end{prop}
\textbf{Proof: } We argue by contradiction. Otherwise there exists a sequence $K_\nu \in C^\infty(M,\mathbb{R})$ converging to $0$ in the $C^1$-topology for which there exist $w_\nu \colon \mathbb{R} \to C^\infty(S^1,V)
\times (T_-,T_+)$ solving the gradient flow equation (\ref{grad}) for $\mathcal{A}_{r_\nu}^{H,K_\nu}$
with $r_\nu \in [0,\infty)$ and satisfying the asymptotic conditions (\ref{asymp}), and $s_\nu \in \mathbb{R}$ with the property that
$$w_\nu(s_\nu) \notin C^\infty(S^1,W) \times (S_-,S_+).$$
We consider the shifted gradient flow lines
$$(v_\nu,\tau_\nu)(s):=w_\nu(s+s_\nu),\quad s \in \mathbb{R}.$$
In particular, we have
$$v_\nu(0) \notin C^\infty(S^1,W)$$
or 
\begin{equation}\label{cont1}
\tau_\nu(0) \in (T-,S-] \cup [S_+,T_+).
\end{equation}
In the first case there exists $t_\nu \in S^1$ such that
\begin{equation}\label{cont2}
v_\nu(0,t_\nu) \in V \setminus W.
\end{equation}
From (\ref{grad2}) we see that $(v_\nu,\tau_\nu)$ solves the problem
\begin{equation}\label{grad3}
\left\{\begin{array}{c}
\partial_s v_\nu(s,t)+J\big(v_\nu(s,t)\big)\Big(\partial_t v_\nu(s,t)-\tau_\nu(s) X_{H+\beta_{r_\nu}(s-s_\nu)K_\nu}\big(v_\nu(s,t)\big)\Big)=0\\
\partial_s \tau_\nu(s)+\int_0^1 \big(H+\beta_{r_\nu}(s-s_\nu)K_\nu\big)\big(v_\nu(s,t)\big)dt=0.
\end{array}\right.
\end{equation}
Using the formulas in the proof of Lemma~\ref{actest}, we can estimate the energy of the gradient flow lines
$w_\nu$ using their asymptotic conditions (\ref{asymp}) as
\begin{eqnarray}\label{energy}
E(w_\nu)&=&\int_{-\infty}^\infty\big|\big\|\partial_s w_\nu(s)\big|\big|^2 ds\\ \nonumber
&=&\int_{-\infty}^\infty\big|\big| \nabla \mathcal{A}_{r_\nu}^{H,K_\nu}\big(w_\nu(s),s\big)\big|\big|^2ds\\ \nonumber
&=&\int_{-\infty}^\infty \partial_s \mathcal{A}_{r_\nu}^{H,K_\nu}\big(w_\nu(s),s\big)ds\\ \nonumber
&\leq&T_+\big(\max K_\nu^+|_{\overline{V}}\,+\max K_\nu^-|_{\overline{V}}\big).
\end{eqnarray}
Therefore the energy of $w_\nu$ converges to zero. Since the energy is invariant under time-shift the same
is true as well for all time-shifts of $w_\nu$. In particular, the energy is uniformly bounded. The first equation
in (\ref{grad3}) states that $v_\nu$ satisfies a perturbed Cauchy-Riemann equation. Since $\tau_\nu(s)$ is uniformly bounded by $T_+$ we see that the pertubation satisfies a uniform $C^0$-bound on the compact set
$\overline{V}$. Since the symplectic form $\omega$ is exact, there is no bubbling and therefore the sequence
$v_\nu$ has a $C^1_{\mathrm{loc}}$-convergent subsequence. In view of the second equation in (\ref{grad3}) the
same holds true for $\tau_\nu$. We denote its limit by $w=(v,\tau)$. Since $K_\nu$ converges to $0$ in the $C^1$-topology we infer that $w$ is a gradient flow line of $\mathcal{A}^H$. Since the energy converges to zero it is
a gradient flow line of no energy. Hence it is a critical point of $\mathcal{A}^H$ and is in particular independent of 
the time-variable $s$. In view of $(\ref{asymp})$ the gradient flow lines $w_\nu$ converge asymptotically to critical points in $\mathcal{C}$. Since the energy converges to zero
in the limit, no breaking can occur and we conclude that $w$ belongs actually to $\mathcal{C}$. However, by
(\ref{cont1}) or (\ref{cont2}) we conclude that 
$$\tau \in [T_-,S_-] \cup [S_+,T_+]$$
or there exists $t \in S^1$ such that
$$v(t) \in \overline{V} \setminus W.$$
This contradicts the fact that $w$ belongs to $\mathcal{C}$. This contradiction proves the proposition. \hfill $\square$
\\
\begin{prop}\label{exist}
Let $\mathcal{U}$ be as in Proposition~\ref{neighbour} and $K \in \mathcal{U}$. Then for every $r \in [0,\infty)$
there exists a solution $w$ of the gradient flow equation (\ref{grad}) for $\mathcal{A}_r=\mathcal{A}_r^{H,K}$ satisfying the asymptotic conditions (\ref{asymp}). 
\end{prop}
\textbf{Proof: } We first discuss the case $r=0$. In this case $\mathcal{A}_0^{H,K}$ just coincides with
$\mathcal{A}^H$. In particular, it is independent of time. Therefore the action is strictly decreasing
along gradient flow lines unless they are constant. In the latter case they have to be critical points of
$\mathcal{A}^H$. In view of the asymptotic conditions the action cannot decrease and therefore the moduli space of solutions precisely coincides with $\mathcal{C}$. 
\\ 
We now fix $w_- \in \mathcal{C}$ and consider the moduli problem
$\mathcal{M}$ consisting of pairs $(\rho,w)$ where $\rho \in [0,r]$ and $w \colon \mathbb{R} \to C^\infty(S^1,V)
\times (T_-,T_+)$ such that
$$\partial_s w(s)+\nabla \mathcal{A}_\rho(w(s))=0,\,\,s \in \mathbb{R}, \quad 
\lim_{s \to -\infty} w(s)=w_-, \quad \lim_{s \to \infty}w(s) \in \mathcal{C}.$$
Note that $(0,w_-)$ with $w_-$ interpreted as a constant gradient flow line belongs to the moduli space $\mathcal{M}$. This is the only member of $\mathcal{M}$ with $\rho=0$. Moreover, it is {nondegenerate} since $\mathcal{C}$ is 
Morse-Bott. 
\\ 
We next show that the moduli space $\mathcal{M}$ is compact. Suppose that $(\rho_\nu,w_\nu)$ is a sequence in $\mathcal{M}$. It follows from Proposition \ref{neighbour}
that 
$$w_\nu(s) \in C^\infty(S^1,W) \times (S_-,S_+)$$ 
for every $s \in \mathbb{R}$ and every $\nu \in \mathbb{N}$.
In particular, since the only periodic orbits of $H$ of period less than or equal to $T_+$ contained in
$\overline{W}$ are the one belonging to $\mathcal{C}$, the sequence of gradient flow lines cannot break at the ends.
Moreover, there is no bubbling since $\omega$ is exact and therefore the sequence has a subsequence which converges
to a limit in $\mathcal{M}$. This shows that the moduli space is compact.  
\\ 
We next assume by contradiction that there is no gradient flow line $w$ such that $(r,w) \in \mathcal{M}$. The moduli space can be interpreted as the zero set of a Fredholm section of index zero from a Hilbert manifold into a Hilbert bundle (See e.g. \cite{FrauenfelderNicolls} for precise notions and formulations). Up to a slight perturbation of the section we can assume that the intersection of the Fredholm section with the zero section is transverse, and hence the moduli space is a {one dimensional} manifold
with boundary where the boundary points are the members $(\rho,w)$ for $\rho=0$ and $\rho=r$. Moreover, for that purpose we do not need to perturb it at the boundary, since the only member for $\rho=0$ is $(0,w_-)$ which is
already transversal. Therefore with our assumption that there is no member for $\rho=r$ we get a compact one-dimensional
manifold with precisely one boundary point, which however does not exist. Consequently there has to exist a solution for $r$. This proves the Proposition. \hfill $\square$
\\ 
\textbf{Proof of Theorem\,A: } We choose $\mathcal{U}$ as in Proposition~\ref{neighbour} and let $K \in \mathcal{U}$. By Proposition~\ref{exist} for every $r \in [0,\infty)$ a solution $w_r$ of the gradient flow equation (\ref{grad}) for $\mathcal{A}_r=\mathcal{A}_r^{H,K}$. Arguing as in the proof of Proposition~\ref{exist} we conclude
that there exists a sequence $r_\nu$ going to infinity such that $w_{r_\nu}$ converges to a gradient flow line
$w_\infty \colon \mathbb{R} \to C^\infty(S^1,V) \times (T_-,T_+)$ of $\mathcal{A}^{H+K}$. Note that in the limit
the action functional $\mathcal{A}^{H+K}$ does not depend on time anymore. Lemma~\ref{actest} tells us that all
gradient flow lines satisfy a uniform action estimate for all time instants,  which then continue to hold for the limit. Therefore we have
for every $\sigma \in \mathbb{R}$
\begin{eqnarray*}
& &\mathcal{A}^H(\mathcal{C})-T_+\big(\max K^+|_{\overline{V}}\,+\max K^-|_{\overline{V}}\big)\\
&\leq& \mathcal{A}^{H+K}\big(w_\infty(\sigma)\big)\\
&\leq& \mathcal{A}^H(\mathcal{C})+T_+\big(\max K^+|_{\overline{V}}\,+\max K^-|_{\overline{V}}\big).
\end{eqnarray*}
In particular, the energy of $w_\infty$ is bounded. Therefore there exists a sequence $s_\nu$ going
to infinity such that $w_\infty(s_\nu)$ converges to a critical point of $\mathcal{A}^{H+K}$. This critical point
then has to satisfy the above action estimate as well. The Theorem is proved. \hfill $\square$

\newpage
\section{Action-Angle Coordinates of the Kepler Problem}\label{Appendix B}
The $d$-dimensional Kepler Problem has a ``hidden'' $SO(d+1)$-symmetry and is super-integrable for $d \ge 2$.  A way to see this is via Moser regularization. In this Appendix we remark on the role of the symmetric group of the Kepler problem in the construction of its action-angle coordinates. Concretely we consider the Kepler problem in dimension 2 with negative energy and we shall design a family of Delaunay-like coordinates for them, by simply considering the symmetric group action on the orbit space with fixed semimajor axis, with the same idea as has been used in Section \ref{Sec: 2}. The same idea extend to Kepler problem in dimension 3 and systems of decoupled Kepler problems in dimension 2 and 3.  We hope to address these in another work. 

Recall that the Delaunay variables of the planar Kepler problem with Hamiltonian
$$H:=\dfrac{\|p\|^{2}}{2}-\dfrac{1}{\|q\|}$$
is the set of canonical coordinates $(L, l, G, g)$ such that 
$$L=\sqrt{a}, \,\,G=\sqrt{a} \sqrt{1-e^{2}},$$
in which $a, e$ are the semimajor axis and the eccentricity of the elliptic Keplerian orbit respectively. The argument of the pericenter $g$ is the angle from the first coordinate axis to the pericenter direction of the orbit and the angle $l$ is the mean anomaly. These variables are well-defined as long as $e \in (0, 1)$. When $e=1$, the orbit is collisional and is non-compact. To properly treat such a situation a regularization procedure has to be involved which we do not address in this Appendix. When $e=0$, the orbit is circular and there are no distinguished pericenter direction, thus the angle $g$ and consequently the angle $l$ are not defined. 

As in Section \ref{Sec: 2} we see that the variable $L=\sqrt{a}$ is well-defined and generates a Hamiltonian circle action for any orbit with $e \in [0, 1)$. The symplectically reduced space $\Omega^{L}$ with respect to this circle action can then be realized as two open hemispheres (according to orientations of the Keplerian motions in the plane) in the sphere $S^{2}$ equipped with an $SO(3)$-invariant symplectic form $\omega_{L}$ which satisfies $\omega_{L}=d G \wedge d g$ in the two open hemispheres. The equator separating the two hemispheres consists of rectilinear orbits and the poles corresponds to circular orbits. 

A realization of this, as in \cite{Albouy}, is to take the unit sphere $S^{2} \subset \R^{3}$ and as orbital plane the horizontal plane $\R^{2} \times \{0\} \subset \R^{3}$. Any point $(x_{1}, x_{2}, x_{3}) \in S^{2}$ determines a bounded Keplerian orbit with semimajor axis $1$, with eccentricity vector $(x_{1}, x_{2})$ and angular momentum $x_{3}$. By rescaling it then gives a bounded Keplerian orbit with fixed semimajor axis $a$ in the horizontal plane. By rotation this argument works for any pre-assigned orbital plane.

With this realization we see that the variables $(G, g)$ are just the symplectic cylindral coordinates of the sphere with respect to the vertical axis of symmetry. 

Now our remark is simply that we may as well pass to any other symplectic cylindral coordinates of the sphere with respect to any axis. In particular we may tilt the axis of symmetry slightly so that the resulting symplectic cylindral coordinates are well-defined for the poles of the sphere, corresponding to circular motions. Denote by $(\tilde{G}, \tilde{g})$ any of such coordinates and the set of variables $(L, l, \tilde{G}, \tilde{g})$ is now canonical as long as they are well-defined. For circular motions however, the angle $l$ is still not well-defined. To remedy this we take the plane in $\R^{3}$ orthogonal to the tilted axis of symmetry. The poles which correspond to circular motions in the horizontal plane now determines an elliptic motion with eccentricity $e \in (0, 1)$ in the tilted plane, to which associated a well-defined set of Delaunay variables $(L, \tilde{l}, \tilde{G}, \tilde{g})$, which we may use as coordinates for the corresponding motions in the horizontal plane.

We now determine the angle $l'$ for near-circular motions.  In a neighborhood of circular motions, a trick of Poincar\'e applies, which leads to a set of symplectic coordinates, see e.g. \cite{Chenciner}. We consider the direct circular motion, which may be characterized by the condition $L=G$. The retrograde circular motion ($L=-G$) can be treated similarly. We consider the angle $\lambda=l+g$ which is well-defined for direct circular motions and motions close by, and further pass to the set of canonical variables
$$(L, \lambda, \xi, \eta)$$
in which $(\xi+ i \eta)=\sqrt{2(L-G)}\exp(-i g)$. These variables are well-defined for direct circular (corresponds to $\xi+i \eta=0$) orbit and orbits close-by. Moreover they are well-defined only except for the retrograde circular motion, which allows us to use to compare different variables. We now pass to Poincar\'e coordinates, $(L, \lambda, \xi, \eta)$ for motions in the horizontal plane and its auxiliary referent motion in the tilted plane $(L, \tilde{\lambda}, \tilde{\xi}, \tilde{\eta})$. The mapping $(L, \lambda, \xi, \eta) \mapsto (L, \tilde{\lambda}=\lambda, \tilde{\xi}, \tilde{\eta})$ is therefore a local canonical transformation for near circular orbits, in which
the angle $\tilde{l}$ is determined by the relationship
$$\tilde{l}+\tilde{g}=l+g=\lambda.$$

By symmetric considerations, we may therefore take $(L, \tilde{l}, \tilde{G}, \tilde{g})$ as  action-angle variables for the corresponding Keplerian motions in the horizontal plane. In particular, they are well-defined near direct circular motions in the horizontal plane. 

\newpage

\section{Restricted Three-Body Problems as Forced Kepler Problems}\label{Appendix C}
A smooth 1-periodic orbit of two massive particles  
$$S^{1} \to \R^{d} \times \R^{d}, t \mapsto (X_{1}(t), X_{2}(t))=(X_{1, 1} (t), \cdots, X_{1, d} (t), X_{2, 1} (t), \cdots, X_{2, d} (t))$$ 
with positive masses $(M_{1}, M_{2})$ and such that $X_{1}(t) \neq X_{2}(t), \forall t \in S_{1}$ determines a restricted three-body problem with Hamiltonian in extended phase space
$$H:=\tau+\dfrac{\|p\|^{2}}{2}-\dfrac{M_{1}}{\|q-X_{1}(t)\|}--\dfrac{M_{2}}{\|q-X_{2}(t)\|},$$
in which we have written $(p, q) =(p_{1}, \cdots, p_{d}, q_{1}, \cdots, q_{d}) \in \R^{d} \times \R^{d}$
with the tautological one form on the cotangent bundle
$$\sum_{i=1}^{d} p_{i} d q_{i} + \tau d t. $$

Setting 
$$p=\tilde{p}+X_{1}'(t), q=\tilde{q}+X_{1}(t)$$
and 
$$\tau=\tilde{\tau}-\sum_{i=1}^{d} \tilde{p}_{i} X_{1, i}'(t)-\|X_{1}'(t)\|^{2} + \sum_{1=1}^{d} \tilde{q}_{i} X_{1,i}^{''}(t), $$ 
we see that 
$$\sum_{i=1}^{d} p_{i} d q_{i} + \tau d t =\sum_{i=1}^{d} p_{i} d \tilde{q}_{i} + \tilde{\tau} d t + d \Bigl(\sum_{i=1}^{d} \tilde{q}_{i} X_{1, i}'(t)\Bigr). $$
Therefore $(\tilde{p}, \tilde{q}, \tilde{\tau}, t)$ form a set of canonical variables on the extended phase space in which the Hamiltonian takes the form
$$H:=\tilde{\tau}+\dfrac{\|\tilde{p}\|^{2}}{2}-\dfrac{M_{1}}{\|\tilde{q}\|}--\dfrac{M_{2}}{\|\tilde{q}+X_{1}(t)-X_{2}(t)\|}-\dfrac{1}{2} \|X_{1}^{'} (t)\|^{2}+\sum_{i=1}^{d} \tilde{q}_{i} X^{''}_{1,i}(t),$$
which assumes the form of \eqref{eq: eq motion perturbed Kepler varepsilon=1} after $M_{1}$ being further normalized to $1$.

Assuming in addition that $(X_{1}(t), X_{2}(t))$ solves a two-body problem (so that each $X_{1}(t)$ moves on a circular or elliptic Keplerian orbit with eccentricity $e \in [0, 1)$ we then obtain the circular and elliptic restricted three-body problem in $\R^{d}$, to which Thm. \ref{thm: main} can be applied.


\begin{thebibliography}{99}

\bibitem{Albouy} A.\,Albouy, Lectures on the two-body problem, in \emph{Classical and Celestial Mechanics: The Recife Lectures}, Princeton University Press, (2002).

\bibitem{Ambrosetti-Coti Zelati} A. Ambrosetti, V. Coti Zelati, \emph{Periodic solutions of singular Lagrangian systems}, Birkh\"auser, (1993).

\bibitem{Ambrosetti-Struwe} A. Ambrosetti, M. Struwe, Periodic motions of conservative systems with singular potentials, \emph{NoDEA} \textbf{1}:179-202, (1994).

\bibitem{AudinDamian} M. Audin, M. Damian, \emph{Morse Theory and Floer Homology}, Springer (2014), Translation by R. Ern\'e from French (EDP Sciences, 2010).

\bibitem{Bahri-Rabinowitz} A. Bahri, P. Rabinowitz, A minimax method for a class of Hamiltonian systems with singular potentials, \emph{J. Funct. Anal.} \textbf{82}: 412-428, (1989).

\bibitem{Verzini-Ortega-Barutello} V.\,Barutello, R.\,Ortega, G.\,Verzini, Regularized variational principles for the perturbed Kepler problem, arXiv:2003.09383, (2020).

\bibitem{Boscaggin-D'Ambrosio-Papini} A.\,Boscaggin, W.\,D'Ambrosio, D.\,Papini, Periodic solutions to a forced Kepler problem in the plane, \emph{Proc. Amer. Math. Soc.}, \textbf{148}:301-314, (2020). 

\bibitem{Boscaggin--Ortega--Zhao} A.\,Boscaggin, R.\,Ortega, L.\,Zhao, Periodic solutions and regularization of a Kepler problem with time-dependent perturbation, \emph{Trans. Amer. Math. Soc.}, \textbf{372}:677-703, (2019). 

\bibitem{Chenciner} A. Chenciner, Int\'egration du probl\`eme de Kepler par la m\'ethode de Hamilton-Jacobi : coordonn\'ees "action-angles" de Delaunay, \emph{Notes Scientifiques et Techniques du Bureau des Longitudes}, S026, (1989).

\bibitem{Cieliebak-Frauenfelder} K. Cieliebak, U. Frauenfelder, A Floer homology for exact contact embeddings, \emph{Pac. J. Math.} \textbf{239}:251-316, (2009).

\bibitem{Frauenfelder-Koert} U.\,Frauenfelder, O.\,van\,Koert, \emph{The Restricted Three-Body Problem and Holomorphic Curves}, Pathways in Mathematics,  Birkh\"auser, (2018).

\bibitem{FrauenfelderNicolls} U. Frauenfelder, R. Nicolls, The moduli space of gradient flow lines and Morse homology, arXiv:2005.10799, (2020).

\bibitem{GinzburgHeinHryniewiczMacarini}, V. L. Ginzburg, D. Hein, U. Hryniewicz, L. Macarini, Closed Reeb orbits on the sphere and symplectically degenerate maxima, \emph{Acta Math. Viet.} \textbf{38}: 55-78, (2013).

\bibitem{GinzburgGurel} V. L. Ginzburg, B. Z. G\"urel, Local Floer homology and the action gap, \emph{J. Sympl. Geom.} \textbf{8}(3):323-357, (2010).

\bibitem{Goursat} E. Goursat, Les transformations isogonales en M\'ecanique. \emph{Comptes Rendus des Séances de l'Académie des Sciences}, \emph{108}:446-448, (1889).

\bibitem{HryniewiczMacarini} U. Hryniewicz, L. Macarini, Local contact homology and applications, \emph{J. Top. Anal.}, \textbf{7}(2):167-238, (2015).

\bibitem{Kummer} M. Kummer, On the regularization of the Kepler problem, \emph{Comm. Math. Phys.}  \textbf{84}(1):133-152, (1982).

\bibitem{Levi-Civita 1904} T.\,Levi-Civita, \emph{Sur la r\'esolution qualitative du probl\`me restreint des trois corps}, ICM report, (1904).

\bibitem{Levi-Civita} T. Levi-Civita, Sur la r\'egularisation du probl\`eme des trois corps. \emph{Acta Math.} \textbf{42}: 99-144, (1920). 

\bibitem{McDuff-Salamon} D. McDuff, D. Salamon, \emph{J-holomorphic Curves and Symplectic Topology, Second Edition}, American Mathematical Society, (2012).

\bibitem{McLean} M. McLean, Local Floer homology and infinitely many simple Reeb orbits, \emph{Algebr. Geom. Topol.}, \textbf{12}(4):1901-1923, (2012).

\bibitem{Moser} J. Moser, Regularization of Kepler's Problem and the Averaging Method on a Manifold, \emph{Comm. Pure Appl. Math.}, \textbf{23}:609-636, (1970).

\bibitem{KimKimKwon} J. Kim, S. Kim, M. Kwon, Bifurcations of symmetric periodic orbits via Floer homology, \emph{Cal. Var. Part. Diff. Eq. } \textbf{59}: Article 101, (2020).


\bibitem{OrtegaZhao}R. Ortega, L. Zhao, Generalized periodic orbits in some restricted three-body problems, arXiv:2007.01619, (2020).

\bibitem{Salamon} D. Salamon, \emph{Lecture on Floer Homology}, \url{https://people.math.ethz.ch/~salamon/PREPRINTS/floer.pdf}, (1997).

\bibitem{Weinstein} A.\,Weinstein, {Normal modes for nonlinear Hamiltonian systems}, \emph{Invent. Math.} \textbf{20}:47-57, (1973).

\bibitem{Zhao} L.\,Zhao, \emph{On some collisional solutions of the rectilinear periodically forced Kepler problem}, \emph{Adv. Nonlinear Stud.}, 16(1): 45-49, (2016).

\end{thebibliography}
\end{document}